\documentclass[11pt]{article} 
\usepackage{latexsym} 
\usepackage[all]{xy} 
\usepackage{amsmath} 
\usepackage{amssymb} 
\usepackage{amsfonts} 
 
\newtheorem{thm}{Theorem}[section] 
\newtheorem{prop}[thm]{Proposition} 
\newtheorem{lem}[thm]{Lemma} 
\newtheorem{df}[thm]{Definition} 
 
\newtheorem{conj}[thm]{Conjecture} 
 
\newtheorem{rmk}[thm]{Remark} 
 
\newtheorem{Q}[thm]{Question} 
 
\begin{document}

\title{\textbf{From HAG to DAG: \\
derived moduli stacks
}} 
\bigskip 
\bigskip 
 
\author{\bigskip\\  
Bertrand To\"en \\  
\small{Laboratoire J. A. Dieudonn\'{e}}\\  
\small{UMR CNRS 6621} \\  
\small{Universit\'{e} de Nice Sophia-Antipolis}\\  
\small{France}\\  
\bigskip \and  
\bigskip \\  
Gabriele Vezzosi \\  
\small{Dipartimento di Matematica}\\  
\small{Universit\`a di Bologna}\\  
\small{Italy}\\  
\bigskip}  

\date{October 2002} 
 
\maketitle 

\begin{abstract} 
These are expanded notes of some talks given during the fall $2002$, about \textit{homotopical algebraic
geometry} with special emphasis on its applications to \textit{derived algebraic geometry}
and \textit{derived deformation theory}. 

We use the general framework developed in \cite{hag1}, and in particular the notions of 
model topology, model sites and stacks over them, in order to define
various \textit{derived moduli functors} and study their geometric properties. 
We start by defining the model category of $D$-stacks, with respect to an extension of 
the \'etale topology to the category of commutative differential graded algebras, and we show
that its homotopy category contains interesting objects, such as schemes, algebraic 
stacks, higher algebraic stacks, dg-schemes, etc. We define the  notion of \textit{geometric
$D$-stacks} and present some related geometric constructions ($\mathcal{O}$-modules,
perfect complexes, $K$-theory, derived tangent stacks, cotangent complexes, various notions of 
smoothness, etc.). Finally, we define and study the derived moduli problems classifying
\textit{local systems} on a topological space, \textit{vector bundles} on a smooth projective variety,
and $A_{\infty}$-\textit{categorical structures}. We state geometricity and smoothness 
results for these examples.

The proofs of the results presented in this paper will be mainly given in \cite{hag2}.
\end{abstract}
 
 
\tableofcontents 
 
\newpage 
 
\begin{section}{Introduction}

These are expanded notes of some talks given during Fall $2002$ about \textit{homotopical algebraic
geometry} with special emphasis on its applications to \textit{derived algebraic geometry}
and \textit{derived deformation theory}. We have omitted proofs that will appear 
mainly in \cite{hag2}. The main purpose of this work is to present in a somehow informal
way the category of $D$-stacks and to give some examples of \textit{derived moduli spaces} as $D$-stacks.\\

We would like to thank the organizers of the conferences
\textit{Axiomatic and enriched homotopy theory} (Cambridge, September 2002) and \textit{Intersection theory and moduli}
(Trieste September 2002), in which some of the material in the present note has been reported. 
We would also like to thank S. M$\mathrm{\ddot{u}}$ller-Stach, H. Esnault and E. Viehweg
for inviting us to lecture on our work at a DFG-Schwerpunkt ``Globale Methoden in der Komplexen Geometrie'' in Essen and the  Research in Pairs program at the 
Matematisches Forschunginstitut Oberwolfach for providing the excellent working conditions in which this paper was written.

\begin{subsection}{What's HAG ?}

\textit{ \textbf{Homotopical Algebraic Geometry} (or \textbf{HAlgebraic Geometry}, or simply \textbf{HAG})
was conceived as a framework to talk about schemes in a context where affine objects are in 
one-to-one correspondence with commutative monoid-like objects in a base symmetric monoidal model category.} \\

This general definition might seem somewhat obscure, so we'd rather mention the most
important examples of base symmetric monoidal model category, and the corresponding notion of
commutative monoid-like objects. In each of the following situations, HAG will provide a context in which one
can \textit{do algebraic geometry} (and in particular, talk about schemes, algebraic spaces, stacks \dots), hence giving rise to various \textit{geometries}.

\begin{enumerate}

\item The model category $Ab$ of abelian groups (with its trivial model structure) and the tensor product
of abelian groups. Commutative monoid objects are commutative rings. The corresponding geometry
is the usual, Grothendieck-style algebraic geometry.

\item The model category $Mod(\mathcal{O})$ of $\mathcal{O}$-modules over some ringed site (with 
the trivial model structure) and the tensor product of $\mathcal{O}$-modules. 
Commutative monoid objects are sheaves of commutative $\mathcal{O}$-algebras. The corresponding geometry
is called \textit{relative algebraic geometry}, and was  
introduced and studied in \cite{hak,del}.

\item The model category $C(k)$ of complexes over some ring k and the tensor product of complexes (see \cite[\S 2.3]{ho}). 
Commutative monoid-like objects are commutative $E_{\infty}$-algebras over $k$ (\cite{km}). The corresponding geometry
is the so-called \textit{derived algebraic geometry} that we are going to discuss in details in this paper, and 
for which one possible avatar is the theory of dg-schemes and dg-stacks of \cite{ck1,ck2}.

\item The model category $Sp$ of 
symmetric spectra and the smash product (see \cite{hss}), or equivalently the category of
 $\mathbb{S}$-modules (see \cite{ekmm}). Commutative monoid-like objects are
$E_{\infty}$-ring spectra, or commutative $\mathbb{S}$-algebras. 
We call the corresponding geometry \textit{brave new algebraic
geometry}, quoting the expression \textit{brave new algebra} introduced
by F. Waldhausen (for more details on the subject, see e.g. \cite{vo}).

\item The model category $Cat$ of categories and the
direct product (see, e.g. \cite{jo}). Commutative monoid-like objects are symmetric monoidal categories.
The corresponding geometry does not have yet a precise name, but could be called \textit{2-algebraic geometry}, 
since vector bundles in this setting will include both the notion of $2$-vector spaces (see \cite{kav}) and its generalization to 
\textit{$2$-vector bundles}.

\end{enumerate}

For the general framework, we refer the reader to \cite{hag1,hag2}. The purpose of the present 
note is to present one possible incarnation of 
HAG through a concrete application to 
\textit{derived algebraic geometry} (or ``DAG'' for short).

\end{subsection}

\begin{subsection}{What's DAG ?}

Of course, the answer we give below is
our own limited understanding of the subject. \\

As far as we know, the foundational ideas of \textit{derived algebraic geometry} (whose infinitesimal theory is also 
referred to as \textit{derived deformation theory}, or ``DTT'' for short)  were introduced by
P. Deligne, V. Drinfel'd and M. Kontsevich, for the purpose of studying the
so-called \textit{derived moduli spaces}. One of the main observation was that certain 
moduli spaces were very \textit{singular} and not of the \textit{expected dimension}, and according to the
general philosophy this was considered as somehow unnatural (see the \textit{hidden smoothness
philosophy} presented in \cite{ko}). It was therefore expected that these
moduli spaces are only \textit{truncations} of some richer geometric objects, 
called the \textit{derived moduli spaces}, containing important additional structures making them
\textit{smooth and of the expected dimension}. In order to illustrate these general ideas, we
present here the fundamental example of the moduli stack of vector bundles (see the introductions of
\cite{ck1,ck2,ka} for more motivating examples as well as philosophical remarks). \\

Let $C$ be a smooth projective curve (say over $\mathbb{C}$), 
and let us consider the moduli
stack $\underline{Vect}_{n}(C)$ of rank $n$ vector bundles on $C$ (here $\underline{Vect}_{n}(C)$ classifies all vector bundles
on $C$, not only the semi-stable or stable ones). The stack $\underline{Vect}_{n}(C)$ is known
to be an algebraic stack (in the sense of Artin). Furthermore, if $E \in \underline{Vect}_{n}(C)(\mathbb{C})$
is a vector bundle on $C$, one can easily compute the \textit{stacky tangent space} of $\underline{Vect}_{n}(C)$
at the point $E$. This \textit{stacky tangent space} is actually a complex of $\mathbb{C}$-vector spaces
concentrated in degrees $[-1,0]$, which is easily seen to be quasi-isomorphic to the complex
$C^{*}(C_{Zar},\underline{End}(E))[1]$ of Zariski cohomology of $C$ with coefficient in the vector
bundle $\underline{End}(E)=E\otimes E^{*}$. Symbolically, 
one writes
$$T_{E}\underline{Vect}(C) \simeq H^{1}(C,\underline{End}(E)) - H^{0}(C,\underline{End}(E)).$$
This implies in particular that the \textit{dimension} of $T_{E}\underline{Vect}(C)$
is independent of the point $E$, and is equal to $n^{2}(g-1)$, where $g$ is the genus of $C$. 
The conlcusion is then that the stack $\underline{Vect}_{n}(C)$ is smooth of dimension $n^{2}(g-1)$.

Let now $S$ be a smooth projective surface, and $\underline{Vect}_{n}(S)$ the moduli stack
of vector bundles on $S$. Once again, $\underline{Vect}_{n}(S)$ is an algebraic stack, and the stacky tangent space 
at a point $E\in \underline{Vect}_{n}(S)(\mathbb{C})$ is easily seen to be given by the same formula
$$T_{E}\underline{Vect}_{n}(S) \simeq H^{1}(S,\underline{End}(E)) - H^{0}(S,\underline{End}(E)).$$
Now, as $H^{2}(S,\underline{End}(E))$ might jump when specializing $E$, the dimension
of $T_{E}\underline{Vect}(S)$, which $h^{1}(S,\underline{End}(E))-h^{0}(S,\underline{End}(E))$,
is not locally constant and therefore the
stack $\underline{Vect}_{n}(S)$ is not smooth anymore.

The main idea of \textit{derived algebraic geometry} is that $\underline{Vect}_{n}(S)$ is only the
truncation of a richer object $\mathbb{R}\underline{Vect}_{n}(S)$, called the \textit{derived
moduli stack of vector bundles on $S$}. This derived moduli stack, whatever it may be, should be such that 
its \textit{tangent space} at a point $E$ is the \textit{whole} complex $C^{*}(S,\underline{End}(E))[1]$, or
in other words, 
$$T_{E}\mathbb{R}\underline{Vect}_{n}(S) \simeq -H^{2}(S,\underline{End}(E)) + 
H^{1}(S,\underline{End}(E)) - H^{0}(S,\underline{End}(E)).$$
The dimension of its tangent space at $E$ is then expected to be $-\chi(S,\underline{End}(E))$, and therefore 
locally constant. Hence, the object $\mathbb{R}\underline{Vect}_{n}(S)$ is expected to be \textit{smooth}. \\

\begin{rmk}\label{r1}
\emph{Another, very similar but probably more striking example 
is given by the moduli stack of stable maps, introduced in \cite{ko}. A consequence of the expected existence of the 
\textit{derived moduli stack of stable maps} is the presence of a  
\textit{virtual structure sheaf} giving rise to a \textit{virtual fundamental class}
(see \cite{bef}). The importance
of such constructions in the context of Gromov-Witten theory shows that the extra information
contained in \textit{derived moduli spaces} is very interesting and definitely geometrically meaningful.} 
\end{rmk}

In the above example of the stack of vector bundles, the tangent space of $\mathbb{R}Vect_{n}(S)$ is expected to be
a complex concentrated in degree $[-1,1]$. More generally, one can get convinced that tangent spaces
of derived moduli (1-)stacks should be complexes concentrated in degree $[-1,\infty[$ (see \cite{ck1}). 
It is therefore pretty clear that in order to make sense of an object such as $\mathbb{R}\underline{Vect}_{n}(S)$, 
schemes and algebraic stacks are not enough, and one should look for a more general 
definition of \textit{spaces}. This leads to the following general question. \\

\bigskip

\textbf{Problem:} \textit{Provide a framework
in which \emph{derived moduli stacks} can actually be constructed. In particular, construct
the derived moduli stack of vector bundles $\mathbb{R}\underline{Vect}(S)$
discussed above.} \\

\bigskip

Several construction of \textit{formal} derived moduli spaces have appeared in the litterature (see for
example \cite{kos,so}), a general framework for \textit{formal DAG} have been 
developed by V. Hinich in \cite{hin2}, and pro-representability questions were investigated by 
Manetti in \cite{man}. So, in a sense, the formal theory has already been 
worked out, and what remains of the problem above is an approach to \textit{global DAG}. \\

A first approach to the global theory was proposed by M. Kapranov and I. Ciocan-Fontanine, and is based 
on the theory of \textit{dg-schemes} or more generally of \textit{dg-stacks} (see \cite{ck1,ck2}). 
A dg-scheme is, roughly speaking, a scheme together with an enrichement of its structural 
sheaf into commutative differential graded algebras. This enriched structural sheaf 
is precisely the datum encoding the \textit{derived} information. 

This approach has been very successful, and
many interesting derived moduli spaces (or stacks) have already been constructed 
as dg-schemes (e.g. the derived version of the Hilbert scheme, of the Quot scheme, 
of the stack of stable maps, and of the stack of local systems on a space have been defined in 
\cite{ka2,ck1,ck2}). However, this approach have encountered two major problems, already identified 
in \cite[0.3]{ck2}.

\begin{enumerate}

\item The definition of dg-schemes and dg-stacks seems too
rigid for certain purposes. By definition, a dg-scheme is a space obtained by 
\textit{gluing commutative differential graded algebras for the Zariski topology}. It seems
however that certain constructions really require a weaker notion of gluing, as for example 
\textit{gluing differential graded algebras up to quasi-isomorphisms}.

\item The notion of dg-schemes is not very well suited with respect to the functorial point of view, 
as representable functors would have to be defined on the derived category of dg-schemes 
(i.e. the category obtained by formally inverting quasi-ismorphisms of dg-schemes), 
which seems difficult to describe and to work with.
As a consequence, the derived moduli spaces constructed in \cite{ka2,ck1,ck2} do not arise as solution to natural
\textit{derived moduli problems}, and are constructed in a rather ad-hoc way. 

\end{enumerate}

The first of these difficulties seems of a technical nature, whereas the second one
seems more fundamental. It seems a direct consequence of 
these two problems that the derived stack of vector bundles still remains  to be constructed in this framework (see \cite{ka} and
\cite[Rem. 4.3.8]{ck1}).  \\

It is the purpose of this note to show how HAG might be applied to provide a framework for DAG in 
which problems $(1)$ and $(2)$ hopefully disappear. We will show in particular how to make  sense
of various \textit{derived moduli functors} whose \textit{representability} can be
proved in many cases.  \\

\end{subsection}

\end{section}

\begin{section}{The model category of $D$-stacks}

In this Section we will present the construction of a \textit{model category of $D$-stacks}. It will be
our derived version of the category of stacks that is commonly used in moduli theory, and all our examples
of derived moduli stacks will be objects of this category. 

The main idea of the construction is the one used in \cite{hag1}, and consists of adopting systematically
the functorial point of view. Schemes, or stacks, are sheaves over the category
of commutative algebras. In the same way, $D$-stacks will be \textit{sheaves-like objects} on
the category of commutative differential graded algebras. This point of view may probably be justified if
one convinces himself that commutative differential graded algebras \textit{have to be} the affine derived moduli
spaces, and that therefore they are the elementary pieces of the theory that one would like to glue to obtain global
geometric objects.  Another, more down to earth, justification would just be to notice that all of the
\textit{wanted derived moduli spaces} we are aware of, have a reasonable model as an object in our category of
$D$-stacks. 

Before starting with the details of the construction, we would like to mention that 
K. Behrend has independently used a similar approach to DAG that uses the $2$-category
of differential graded algebras (see \cite{be}) (his approach is actually the $2$-truncated version
of ours). It is not clear to us that the constructions and results we are going to
present in this work have reasonable analogs in his framework, as they use in an essential way 
higher homotopical informations that are partially lost when using any truncated version. \\

\textbf{Conventions.} For the sake of simplicity, we will work over the field of
complex numbers $\mathbb{C}$. The expression \textit{cdga} will always refer to a  
\textit{non-positively graded commutative differential graded algebra} over $\mathbb{C}$ with 
differential of degree 1. Therefore, a cdga $A$ looks like
$$\xymatrix{ \dots \ar[r] & A^{-n} \ar[r] & A^{-n+1} \ar[r] & \dots \ar[r] & A^{-1} \ar[r] & A^{0}.}$$
The category CDGA of cdga's is endowed with its usual model category structure (see e.g. \cite{hin}), for which 
fibrations (resp. equivalences) are epimorphisms in degree $\leq -1$ (resp. quasi-isomorphisms). 

\begin{subsection}{$D$-Pre-stacks}

We start by defining $D-Aff:=CDGA^{op}$ to be the opposite category of cdga's, and we consider
the category $SPr(D-Aff)$, of simplicial presheaves on $D-Aff$, or equivalently the category of
functors from $CDGA$ to $SSet$. The category 
$SPr(D-Aff)$ is endowed with its objectwise projective model structure in which fibrations and
equivalences are defined objectwise (see \cite[13.10.17]{hi}). 

For any cdga $A \in D-Aff$, we have the presheaf of sets represented by $A$, denoted by
$$\begin{array}{cccc}
h_{A} : & D-Aff^{op} & \longrightarrow & Set \\
 & B & \mapsto & Hom(B,A).
\end{array}$$
Considering a set as a constant simplicial set, we will look at $h_{A}$ as an object in 
$SPr(D-Aff)$.
The construction $A \mapsto h_{A}$ is clearly functorial in $A$, and therefore for any
$u : A \rightarrow A'$ in $D-Aff$, corresponding to a quasi-isomorphism of cdga's, we
get a morphism $u : h_{A} \rightarrow h_{A'}$ in $SPr(D-Aff)$. These morphisms will simply be 
called \textit{quasi-isomorphisms}. 

\begin{df}\label{d1}
The \emph{model category of $D$-pre-stacks} is the left Bousfield localization of
the model category $SPr(D-Aff)$ with respect to the set of morphisms $\{u : h_{A} \rightarrow h_{A'}\}$,
where $u$ varies in the set of all quasi-isomorphisms.
It will be denoted by $D-Aff^{\wedge}$.
\end{df}

\begin{rmk}
\begin{enumerate}
\item \emph{The careful reader might object that the category $D-Aff$ and the set of all quasi-isomorphisms are not small, and therefore
that Definition \ref{d1} does not make sense. If this happens (and only then), take two universes
$\mathbb{U} \in \mathbb{V}$, define $CDGA$ as the category of $\mathbb{U}$-small cdga's and $SPr(D-Aff)$ as the
category of functors from $CDGA$ to the category of $\mathbb{V}$-small simplicial sets. 
Definition \ref{d1} will now make sense.}
\item \emph{In \cite{hag1}, the model category $D-Aff^{\wedge}$ was denoted by $(D-Aff,W)^{\wedge}$, where $W$
is the subcategory of quasi-isomorphisms.}
\end{enumerate}
\end{rmk}

By general properties of left Bousfield localization (see \cite{hi}), the fibrant objects in 
$D-Aff^{\wedge}$ are the functors $F : CDGA \longrightarrow SSet$ satisfying the following two conditions

\begin{enumerate}
\item For any $A \in CDGA$, the simplicial set $F(A)$ is fibrant.

\item  For any quasi-isomorphism $u : A \longrightarrow B$ in CDGA, the induced morphism 
$F(u) : F(A) \longrightarrow F(B)$ is a weak equivalence of simplicial sets.

\end{enumerate}

From this description, we conclude in particular, that the homotopy category $\mathrm{Ho}(D-Aff^{\wedge})$
is naturally equivalent to the full sub-category of $\mathrm{Ho}(SPr(D-Aff))$ consisting of
functors $F : CDGA \longrightarrow SSet$ sending quasi-isomorphisms to weak equivalences. We will use
implicitely this description, and we will always consider $\mathrm{Ho}(D-Aff^{\wedge})$ as embedded in 
$\mathrm{Ho}(SPr(D-Aff))$. 

\begin{df}\label{d1'}
Objects of $D-Aff^{\wedge}$ satisfying condition $(2)$ above (i.e. sending 
quasi-isomorphisms to weak equivalences) will be
called $D$\emph{-pre-stacks}.
\end{df}

\end{subsection}

\begin{subsection}{$D$-Stacks}

Now that we have constructed the model category of $D$-pre-stacks we will 
introduce some kind of \textit{\'etale topology} on the category $D-Aff$. This will
allow us to talk about a corresponding notion of \'etale \textit{local equivalences}
in $D-Aff^{\wedge}$, and to define the model category of $D$-stacks by including
the \textit{local-to-global principle} into the model structure. \\

We learned the following notion of formally \'etale morphism of cdga's from K. Behrend.

\begin{df}\label{d2}
A morphism $A \longrightarrow B$ in $CDGA$ is called \emph{formally \'etale}
if it satisfies the following two conditions.
\begin{enumerate}
\item The induced morphism $H^{0}(A) \longrightarrow H^{0}(B)$ is a formally \'etale morphism
of commutative algebras.

\item For any $n<0$, the natural morphism of $H^{0}(B)$-modules
$$H^{n}(A)\otimes_{H^{0}(A)}H^{0}(B) \longrightarrow H^{n}(B)$$
is an isomorphism.

\end{enumerate}
\end{df} 

\begin{rmk}
\emph{Its seems that a morphism $A \longrightarrow B$ 
of cdga's is formally \'etale in the sense of Definition \ref{d2} if and only if
the relative cotangent complex $\mathbb{L}\Omega^{1}_{B/A}$ (e.g. in the sense of \cite{hin}) is acyclic. This 
justifies the terminology.}
\end{rmk}

From Definition \ref{d2} we now define the notion of \'etale covering families. For this, we recall that 
a morphism of cdga's $A \longrightarrow B$ is said to be \textit{finitely presented} if
$B$ is equivalent to a retract of a finite cell $A$-algebra (see for example \cite{ekmm}). This is also equivalent to say that
for any filtered systems $\{A \longrightarrow C_{i}\}_{i\in I}$, the natural morphism 
$$Colim_{i\in I} \mathrm{Map}_{A/CDGA}(B,C_{i}) \longrightarrow \mathrm{Map}_{A/CDGA}(B,Colim_{i\in I}C_{i})$$
is a weak equivalence (here $\mathrm{Map}_{A/CDGA}$ denotes the mapping spaces, or function complexes, of the
model category $A/CDGA$ of cdga's under $A$, as defined in \cite[\S 5.4]{ho})\footnote{
We warn the reader that if commutative algebras are considered as cdga's concentrated in degree
zero, the notion of finitely presented morphisms of commutative algebras and
the notion of finitely presented morphisms of cdga's are \emph{not the same}. In fact, for a morphism 
of commutative algebras it is stronger to be finitely presented as a morphism of cdga's than
as a morphism of algebras.}.

\begin{df}\label{d3}
A finite family of morphisms of cdga's
$$\{A \longrightarrow B_{i}\}_{i \in I}$$
is called an \emph{\'etale covering} if it satisfies the following three conditions
\begin{enumerate}
\item For any $i \in I$, the morphism $A \longrightarrow B_{i}$ is
finitely presented.

\item For any $i \in I$, the morphism $A \longrightarrow B_{i}$ is
formally \'etale.

\item The induced family of morphisms of affine schemes
$$\{Spec\, H^{0}(B_{i}) \longrightarrow Spec\, H^{0}(A)\}_{i \in I}$$
is an \'etale covering. 

\end{enumerate}
\end{df}

The above definition almost defines a pre-topology on the category $D-Aff$. Indeed, stability and 
composition axioms for a pre-topology are statisfied, but the base change axiom is not. 
In general, the base change of an \'etale covering $\{A \longrightarrow B_{i}\}_{i \in I}$ along 
a morphism of $A \longrightarrow C$ will only be an \'etale covering 
if $A \longrightarrow C$ is a cofibration in CDGA. In other words, for the base change axiom
to be satisfied one needs to replace fibered products by homotopy fibered products in 
$D-Aff$. Therefore, the \'etale covering families
of Definition \ref{d3} do not satisfy the axioms for a pre-topology on $D-Aff$, but rather satisfy 
a \textit{homotopy analog} of them. This is an example of a \textit{model pre-topology}
on the model category $D-Aff$, for which we refer the reader to \cite[\S 4.3]{hag1} where  a precise definition is given. 

In turns out that the data of a model pre-topology on a model category $M$ is more or
less equivalent to the data of a Grothendieck topology on its homotopy category $\mathrm{Ho}(M)$ (see \cite[Prop. 4.3.5]{hag1}). 
In our situation, the \'etale coverings of Definition \ref{d3} induce 
a Grothendieck topology, called the \textit{\'etale topology} on the 
opposite of the homotopy category $\mathrm{Ho}(D-Aff)$ of cdga's. More concretely, a sieve
$S$ over a cdga $A \in \mathrm{Ho}(D-Aff)$ is declared to be a covering sieve if
it contains an \'etale covering family $\{A \longrightarrow B_{i}\}_{i \in I}$. 
The reader will check as an exercise that this defines a topology on $\mathrm{Ho}(D-Aff)$
(hint: one has to use that \'etale covering families are stable by homotopy pull-backs in $D-Aff$, 
or equivalentely by homotopy push-outs in $CDGA$). From now on, we will 
always consider $\mathrm{Ho}(D-Aff)$ as a Grothendieck site for this \'etale topology. \\

For a $D$-pre-stack $F \in D-Aff^{\wedge}$ (recall from Definition \ref{d1'} that this
implies that $F$ sends quasi-isomorphisms to weak equivalences), we define
its presheaf of connected components
$$\begin{array}{cccc}
\pi_{0}^{pr}(F) : & D-Aff^{op} & \longrightarrow & Set \\
& A & \mapsto & \pi_{0}(F(A)).
\end{array}$$
As the object $F$ is a $D$-pre-stack (see \ref{d1'}), the functor $\pi_{0}^{pr}(F)$ 
will factors
 through the homotopy category
$$\begin{array}{cccc}
\pi_{0}^{pr}(F) : & \mathrm{Ho}(D-Aff)^{op} & \longrightarrow & Set \\
& A & \mapsto & \pi_{0}(F(A)).
\end{array}$$
We can consider the sheaf $\pi_{0}(F)$ associated to 
the presheaf $\pi_{0}^{pr}$ in the \'etale topology on $\mathrm{Ho}(D-Aff)$. The sheaf $\pi_{0}(F)$  
is called the $0$-\textit{homotopy sheaf} of the $D$-pre-stack $F$. 
Now, if $F \in D-Aff^{\wedge}$ is any simplicial presheaf, then one can apply the above
construction to one of its fibrant models $RF$. This allows us to define its
$0$-th homotopy sheaf as $\pi_{0}(F):=\pi_{0}(RF)$.

As for the case of simplicial presheaves (see \cite{ja1}), one can also define
higher homotopy sheaves, which are sheaves of groups and abelian groups
on the sites $\mathrm{Ho}(D-Aff/A)$ for various cdga's $A$. Precisely, let $F$ be a
$D$-pre-stacks and $s \in F(A)_{0}$ a point over a cdga $A \in D-Aff$.
We define the $n$-th homotopy group presheaf pointed at $s$ by
$$\begin{array}{cccc}
\pi_{n}^{pr}(F,s) : & D-Aff^{op}/A=A/CDGA & \longrightarrow & Gp \\
& (u : A \rightarrow B) & \mapsto & \pi_{n}(F(B),u^{*}(s)).
\end{array}$$
Again, as $F$ is a $D$-pre-stack, this presheaves descend to the homotopy category
$$\begin{array}{cccc}
\pi_{n}^{pr}(F,s) : & \mathrm{Ho}(D-Aff^{op}/A)=\mathrm{Ho}(A/CDGA) & \longrightarrow & Gp \\
& (u : A \rightarrow B) & \mapsto & \pi_{n}(F(B),u^{*}(s)).
\end{array}$$
The \'etale model pre-topology on $D-Aff$ also induces Grothendieck
topologies on the various homotopy categories $\mathrm{Ho}(A/CDGA)$, and therefore
one can consider the sheaves associated to $\pi_{n}^{pr}(F,s)$. These
sheaves are called the $n$-\textit{th homotopy sheaves} of $F$ pointed at $s$ and 
are denoted by $\pi_{n}(F,s).$ As before, if $F$ is any object in $D-Aff^{\wedge}$, one can 
define $\pi_{n}(F,s):=\pi_{n}(RF,s)$ for $RF$ a fibrant replacement of $F$. \\

The notion of homotopy sheaves defined above gives rise to the following notion of
local equivalences.

\begin{df}\label{d4}
A morphism $f : F \longrightarrow F'$ in $D-Aff^{\wedge}$ is called a \emph{local equivalence}
if it satisfies the following two conditions
\begin{enumerate}
\item The induced morphism of sheaves $\pi_{0}(F) \longrightarrow \pi_{0}(F')$ is 
an isomorphism.

\item For any $A \in D-Aff$, and any point $s \in F(A)$, the induced morphism of sheaves
$\pi_{n}(F,s) \longrightarrow \pi_{n}(F',f(s))$ is an isomorphism.
\end{enumerate}
\end{df}

One of the key results of ``HAG'' is the following theorem. It is a very special case
of the existence theorem \cite[\S 4.6]{hag1}, which extends the existence of the local model structure
on simplicial presheaves (see \cite{ja1}) to the case of model sites. 

\begin{thm}\label{t1}
There exists a model category structure on $D-Aff^{\wedge}$ for which the
equivalences are the local equivalences and the cofibrations are
the cofibrations in the model category $D-Aff^{\wedge}$ of $D$-pre-stacks.

This model category is called the \emph{model category of $D$-stacks} for the \'etale topology, and is denoted
by $D-Aff^{\sim}$. 
\end{thm}

The reason for calling $D-Aff^{\sim}$ the model category of $D$-stacks is the following
proposition. It follows from \cite[4.6.3]{hag1}, which is a generalization to model sites of the main theorem of \cite{dhi}. 

\begin{prop}\label{p1}
An object $F \in D-Aff^{\sim}$ is fibrant if and only if it satisfies the following
three conditions
\begin{enumerate}
\item For any $A \in D-Aff$, the simplicial set $F(A)$ is fibrant.

\item For any quasi-isomorphism of cdga's $A \longrightarrow B$,
the induced morphism $F(A) \longrightarrow F(B)$ is a weak equivalence.

\item For any cdga $A$, and any \'etale hyper-covering in $D-Aff$ (see \cite{hag1} for details)
$A \longrightarrow B_{*}$, the induced morphism
$$F(A) \longrightarrow Holim_{n \in \Delta}F(B_{n})$$
is a weak equivalence. 

\end{enumerate}
\end{prop}

Condition $(3)$ is called the \textit{stack condition for the \'etale topology}. 
Note that a typical \'etale hyper-covering  of cdga's $A \longrightarrow B_{*}$ is given
by the homotopy co-nerve of an \'etale covering morphism $A \longrightarrow B$ 
$$B_{n}:=\underbrace{B\otimes^{\mathbb{L}}_{A}B\otimes^{\mathbb{L}}_{A} \dots 
\otimes^{\mathbb{L}}_{A}B}_{n\; times}.$$
Condition $(3)$ for these kind of hyper-coverings is the most commonly used 
descent condition, but as first shown in \cite{dhi}
requiring descent with respect to \textit{all} \'etale hyper-coverings is necessary for Proposition \ref{p1} to be correct.

\begin{df}\label{d5}
A $D$\emph{-stack} is any object $F \in D-Aff^{\sim}$ satisfying conditions
$(2)$ and $(3)$ of Proposition \ref{p1}. 
By abuse of language, objects in the homotopy category $\mathrm{Ho}(D-Aff^{\sim})$ will 
also be called $D$-stacks.

A \emph{morphism} of $D$-stacks is a morphism in the homotopy category $\mathrm{Ho}(D-Aff^{\sim})$. 
\end{df}

The second part of the definition is justified because the homotopy category
$\mathrm{Ho}(D-Aff^{\sim})$ is naturally equivalent to the full sub-category of
$\mathrm{Ho}(SPr(D-Aff))$ consisting of objects satisfying conditions $(2)$ and $(3)$ 
of Proposition \ref{p1}. 

\end{subsection}

\begin{subsection}{Operations on $D$-stacks}

One of the main consequences of the existence of the model structure on $D-Aff^{\sim}$ is
the possibility to define several standard operations on $D$-stacks, analogous
to the ones used in sheaf theory (limits, colimits, sheaves of morphisms \dots). \\

First of all, the category $D-Aff^{\sim}$ being a category of simplicial presheaves, 
it comes with a natural enrichement over the category of simplicial sets. This makes
$D-Aff^{\sim}$ into a simplicial model category (see \cite[4.2.18]{ho}). In particular, 
one can define in a standard way the \textit{derived simplicial Hom's} (well defined
in the homotopy category $\mathrm{Ho}(SSet)$),
$$\mathbb{R}\underline{Hom}(F,G):=\underline{Hom}(QF,RG),$$
where $Q$ is a cofibrant replacement functor, $R$ is a fibrant replacement functor, and
$\underline{Hom}$ are the simplicial Hom's sets of $D-Aff^{\sim}$. These
derived simplicial Hom's allows one to consider spaces of morphisms between $D$-stacks, in the same
way as one commonly considers groupoids of morphisms between stacks in groupoids (see \cite{lm}).

This simplicial structure also allows one to define \textit{exponentials} by simplicial sets. 
For an object $F \in D-Aff^{\sim}$ and $K \in SSet$, one has a well defined object in $\mathrm{Ho}(D-Aff^{\sim})$ 

$$F^{\mathbb{R}K}:=(RF)^{K}$$
which satisfies the usual adjunction formula
$$\mathbb{R}\underline{Hom}(G,F^{\mathbb{R}K}) \simeq \mathbb{R}\underline{Hom}(K,\mathbb{R}\underline{Hom}(G,F)).$$

The existence of the model structure $D-Aff^{\sim}$ also implies the existence of \textit{homotopy limits} and \textit{homotopy colimits}, 
as defined in \cite[\S 19]{hi}. The existence of these homotopy limits and colimits 
is the analog of the fact that category of sheaves have all kind of limits and colimits. 
We will use in particular homotopy pull-backs i.e. homotopy limits of diagrams $\xymatrix{F & H \ar[r] \ar[l] & G}$, 
that will be denoted by 
$$F\times_{H}^{h}G:=Holim \{ \xymatrix{F & H \ar[r] \ar[l] & G} \}.$$

Finally, one can show that the homotopy category $\mathrm{Ho}(D-Aff^{\sim})$ is \textit{cartesian closed} (see \cite[\S 4.7]{hag1}). 
Therefore, for any two object $F$ and $G$, there exists an object $\mathbb{R}\mathcal{HOM}(F,G) \in \mathrm{Ho}(D-Aff^{\sim})$,  
which is determined by the natural isomorphisms
$$\mathbb{R}\underline{Hom}(F\times G,H) \simeq \mathbb{R}\underline{Hom}(F,\mathbb{R}\mathcal{HOM}(G,H)).$$
We say that $\mathbb{R}\mathcal{HOM}(F,G)$ is the $D$-\textit{stack of morphisms} from $F$ to $G$, analogous 
to the sheaf of morphisms between two sheaves. 

If one looks at these various constructions, one realizes that
$D-Aff^{\sim}$ has all the homotopy analogs of the properties that characterize Grothendieck topoi. To be more precise, C. Rezk has
defined a notion of \textit{homotopy topos} (we rather prefer the expression \textit{model topos}), 
which are model categories behaving homotopically very much like a usual topos.
The standard examples of such homotopy topoi are model categories of simplicial presheaves on some Grothendieck 
site, but not all
of them are of this kind; the model category $D-Aff^{\sim}$ is in fact an example of a model topos which is 
not equivalent to model categries of simplicial presheaves on some site (see \cite[\S 3.8]{hag1} for more
details on the subject). 

\end{subsection}

\end{section}

\begin{section}{First examples of $D$-stacks}

Before going further with the geometric properties of $D$-stacks, we would like to present
some examples. More examples will be given in the Section 5.

\begin{subsection}{Representables}

The very first examples of schemes are affine schemes. In the same way, our first example
of $D$-stacks are \textit{representable $D$-stacks}\footnote{We could as well have called
them \textit{affine $D$-stacks}.}. 

We start by fixing a fibrant resolution functor $\Gamma$ on the model category $CDGA$. 
Recall that this means that for any cdga $B$, $\Gamma(B)$ is a simplicial object
in $CDGA$, together with a natural morphism $B \longrightarrow \Gamma(B)$ that makes it into 
a fibrant replacement for the Reedy model structure on simplicial objects (see \cite[\S 5.2]{ho}). In the present situation, 
one could choose the following standard fibrant resolution functor
$$\begin{array}{cccc}
\Gamma(B) : & \Delta^{op} & \longrightarrow & CDGA \\
 & [n] & \mapsto & \Gamma(B)_{n}:=B\otimes \Omega^{*}_{\Delta^{n}}.
\end{array}$$
Here $\Omega^{*}_{\Delta^{n}}$ is the cdga (exceptionally positively graded) of 
algebraic differential forms on the standard algebraic $n$-simplex. Of course the cdga $B\otimes \Omega^{*}_{\Delta^{n}}$
is not non-positively graded, but one can always take its truncation in order to see it as an object 
in $CDGA$. 

Now, for any cdga $A$, we define a functor
$$\begin{array}{cccc}
Spec\, A : & CDGA & \longrightarrow & SSet \\
 & B & \mapsto & Hom(A,\Gamma(B)),
\end{array}$$
that is considered as an object in $D-Aff^{\sim}$. This construction is clearly functorial in $A$ and gives
rise to a functor
$$Spec : CDGA^{op}=D-Aff \longrightarrow D-Aff^{\sim}.$$
The functor $Spec$ is almost a right Quillen functor: it preserves fibrations, trivial fibrations and
limits, but does not have a left adjoint. However, it has a well defined right derived functor
$$\mathbb{R}Spec\, : \mathrm{Ho}(CDGA)^{op}=\mathrm{Ho}(D-Aff) \longrightarrow \mathrm{Ho}(D-Aff^{\sim}).$$

A fundamental property of this functor is the following lemma.

\begin{lem}\label{l1}
The functor $\mathbb{R}Spec$ is fully faithful. More generally, for two cdga's $A$ and $B$, it induces
a natural equivalence on the mapping spaces
$$\mathbb{R}\underline{Hom}(A,B) \simeq \mathbb{R}\underline{Hom}(\mathbb{R}Spec\, B,\mathbb{R}Spec\, A).$$ 
\end{lem}

The above lemma contains two separated parts. The first part states that $\mathbb{R}Spec\,$ is 
fully faithful when considered to have values in  $\mathrm{Ho}(D-Aff^{\wedge})$ (i.e. when one forgets about the topology). This first part
is a very general result that we call \textit{Yoneda lemma for model categories} (see \cite[\S 4.2]{hag1}).
The second part of the lemma states that for a cofibrant cdga $A$, the
object $Spec(A)$ is a $D$-stack (see Definition \ref{d5}). This is not a 
general fact, and of course depends on the choice of the topology. Another way to express  this
last result is to say that the \textit{\'etale topology is sub-canonical}. 

\begin{df}\label{d6}
A $D$-stack isomorphic in $\mathrm{Ho}(D-Aff^{\sim})$ to some $\mathbb{R}Spec\, A$ is called
a \emph{representable} $D$-stack.
\end{df}

In particular, Lemma \ref{l1} implies that the full subcategory of $\mathrm{Ho}(D-Aff)^{\sim}$ consisting of
representable $D$-stacks is equivalent to the homotopy category of cdga's. 

\end{subsection}

\begin{subsection}{Stacks vs. $D$-stacks}

Our second example of $D$-stacks are simply stacks. In other words, any stack defined
over the category of affine schemes with the \'etale topology gives rise to a $D$-stack. \\

Let $Alg$ be the category of commutative $\mathbb{C}$-algebras, and $Aff=Alg^{op}$ its opposite 
category. Recall that there exists a model category of simplicial presheaves on 
$Aff$ for the \'etale topology (see \cite{ja1}). We will consider its projective version
described in \cite{bl}, and denote it by $Aff^{\sim}$. This model category is called
the \textit{model category of stacks for the \'etale topology}. Its
homotopy category $\mathrm{Ho}(Aff^{\sim})$ contains as full subcategories 
the category of sheaves of sets and the category of stacks in groupoids (see e.g. \cite{lm}). 
More generally, one can show that the full subcategory of $n$-truncated objects
in $\mathrm{Ho}(Aff^{\sim})$ is naturally equivalent to the homotopy category of
stacks in $n$-groupoids (unfortunately there are no references for this last result until now but  
the reader might consult \cite{hol} for the case $n=1$). In particular, $\mathrm{Ho}(Aff^{\sim})$ contains
as a full subcategory the category of schemes, and more generally of Artin stacks. 

There exists an adjunction
$$H^{0} : CDGA \longrightarrow Alg \qquad CDGA \longleftarrow Alg : j,$$
for which $j$ is the full embedding of $Alg$ in $CDGA$ that sends a commutative algebra
to the corresponding cdga concentrated in degree $0$. Furthermore, this adjunction
is a Quillen adjunction when $Alg$ is endowed with its trivial model structure 
(as written above, $j$ is on the right and $H^{0}$ is its left adjoint). This adjunction induces 
various adjunctions between the category of simplicial presheaves
$$j_{!} : Aff^{\sim} \longrightarrow D-Aff^{\sim} \qquad Aff^{\sim} \longleftarrow D-Aff^{\sim} : j^{*}$$
$$j^{*} : D-Aff^{\sim} \longrightarrow Aff^{\sim} \qquad D-Aff^{\sim} \longleftarrow Aff^{\sim} : (H^{0})^{*}$$
One can check that these adjunction are Quillen adjunction (where the functors written on the left are left Quillen).
In particular we conclude that $j^{*}$ is right and left Quillen, and therefore preserves equivalences. 
From this we deduce easily the following important fact.

\begin{lem}\label{l2}
The functor
$$i:=\mathbb{L}j_{!} : \mathrm{Ho}(Aff^{\sim}) \longrightarrow \mathrm{Ho}(D-Aff^{\sim})$$
is fully faithful. 
\end{lem}

The important consequence of the previous lemma is that $\mathrm{Ho}(D-Aff^{\sim})$ contains 
schemes, algebraic stacks \dots, as full sub-categories. \\

\textbf{Warning:} The full embedding $i$ does not commute with homotopy pull-backs, nor with
internal Hom-$D$-stacks. \\

This warning is the real heart of DAG: the category of $D$-stacks contains usual stacks, but 
these are not stable under the standard operations of $D$-stacks. In other words, if one starts with 
some schemes and performs some constructions on these schemes, considered as $D$-stacks, 
the result might not be a scheme anymore. This is the main reason why derived moduli
spaces are not schemes, or stacks in general ! \\

\textbf{Notations.} In order to avoid confusion, a scheme or a stack $X$, when considered as 
a $D$-stack will always be denoted by $i(X)$, or simply by $iX$. \\

The full emdedding $i=\mathbb{L}j_{!}$ has a right adjoint $\mathbb{R}j^{*}=j^{*}$. It will be
denoted by 
$$h^{0}:=j^{*} : \mathrm{Ho}(D-Aff^{\sim}) \longrightarrow \mathrm{Ho}(Aff^{\sim}),$$
and called the \textit{truncation functor}. Note that for any cdga, one has
$$h^{0}(\mathbb{R}Spec\, A)\simeq Spec\, H^{0}(A),$$
which justifies the notation $h^{0}$. Note also that for any $D$-stack $F$, and any 
commutative algebra $A$, one has
$$F(A)\simeq \mathbb{R}\underline{Hom}(iSpec\, A, F) \simeq 
\mathbb{R}\underline{Hom}(Spec\, A,h^{0}(F))\simeq h^{0}(F)(A).$$
This shows that a $D$-stack $F$ and its truncation $h^{0}(F)$ \textit{have the same
points with values in commutative algebras}. Of course, $F$ and $h^{0}(F)$ do not have
the same points with values in cdga's in general, except when $F$ is of the form $iF'$
for some stack $F' \in \mathrm{Ho}(Aff^{\sim})$. \\

\textbf{Terminology.} Points with values in commutative algebras will be called
\textit{classical points}. \\

\noindent We just saw that a $D$-stack $F$ and its truncation $h^{0}(F)$
always have the same classical points. 

Given two stacks $F$ and $G$ in $Aff^{\sim}$, there exists a stack of morphisms 
$\mathbb{R}\mathcal{HOM}(F,G)$, that is the derived internal Hom's of the model category
$Aff^{\sim}$ (see \cite[\S 4.7]{hag1}). 
As remarked above, the two objects $i\mathbb{R}\mathcal{HOM}(F,G)$
and $\mathbb{R}\mathcal{HOM}(iF,iG)$ are different in general. However, one has
$$h^{0}(\mathbb{R}\mathcal{HOM}(iF,iG))\simeq \mathbb{R}\mathcal{HOM}(F,G),$$
showing that $i\mathbb{R}\mathcal{HOM}(F,G)$ and  $\mathbb{R}\mathcal{HOM}(iF,iG)$ have the same classical points. 

\end{subsection}

\begin{subsection}{dg-Schemes}

We have just seen that the homotopy category of $D$-stacks $\mathrm{Ho}(D-Aff^{\sim})$ contains
the categories of schemes and algebraic stacks. We will now relate the notion
of dg-schemes of \cite{ck1,ck2} to $D$-stacks. \\

Recall that a dg-scheme is a pair $(X,\mathcal{A}_{X})$, consisting of a scheme $X$ and
a sheaf of $\mathcal{O}_{X}$-cdga's on $X$ such that $\mathcal{A}_{X}^{0}=\mathcal{O}_{X}$ (however, this last condition
does not seem so crucial). For the sake of simplicity we will assume that $X$ is 
quasi-compact and separated. We can therefore take a finite affine open covering $\mathcal{U}=\{U_{i}\}_{i}$ of 
$X$, and consider
its nerve $N(\mathcal{U})$ (which is a simplicial scheme)
$$\begin{array}{cccc}
N(\mathcal{U}) : & \Delta^{op} & \longrightarrow & \{Schemes\} \\
 & [n] & \mapsto & \coprod_{i_{0},\dots,i_{n}}U_{i_{0},\dots,i_{n}}
\end{array}$$
where, as usual, $U_{i_{0},\dots,i_{n}}=U_{i_{0}}\cap U_{i_{1}}\cap \dots U_{i_{n}}$. Note that 
as $X$ is separated and the covering is finite, $N(\mathcal{U})$ is in fact a simplicial affine scheme. 

For each integer $n$, let $A(n)$ be the cdga of global sections of $\mathcal{A}_{X}$ on the
scheme $N(\mathcal{U})_{n}$. In other words, one has
$$A(n)=\prod_{i_{0},\dots,i_{n}}\mathcal{A}_{X}(U_{i_{0}})\times \mathcal{A}_{X}(U_{i_{1}}) \times
\dots \mathcal{A}_{X}(U_{i_{n}}).$$

The simplicial structure on $N(\mathcal{U})$ makes $[n] \mapsto A(n)$ into a 
co-simplicial diagram of cdga's. By applying levelwise the functor $\mathbb{R}Spec$, we get a
simplicial object $[n] \mapsto \mathbb{R}Spec\, A(n)$ in $D-Aff^{\sim}$. We define
the $D$-stack $\Theta(X,\mathcal{A}_{X}) \in \mathrm{Ho}(D-Aff^{\sim})$ to be the homotopy colimit of this diagram
$$\Theta(X,\mathcal{A}_{X}):=Hocolim_{[n] \in \Delta^{op}}\mathbb{R}Spec\, A(n).$$

One can check, that $(X,\mathcal{A}_{X}) \mapsto \Theta(X,\mathcal{A}_{X})$ defines a functor
$$\Theta : \mathrm{Ho}(dg-Sch) \longrightarrow \mathrm{Ho}(D-Aff^{\sim}),$$
from the homotopy category of (quasi-compact and separated) dg-schemes to the homotopy category of 
$D$-stacks. This functor allows us to consider dg-schemes as $D$-stacks. \\

\textbf{Question:} \textit{Is the functor $\Theta$ fully faithful ?} \\

We do not know the answer to this question, and there are no real reasons for this
 answer to be positive. As already explained in the Introduction, the difficulty comes from the fact
that the homotopy category of dg-schemes seems quite difficult to describe. 
In a way, it might not be so important to know the answer to the above question, 
as until now morphisms in the homotopy category of dg-schemes have never been taken into account
seriously, and only the objects of the category $\mathrm{Ho}(dg-Sch)$ have been shown to be
relevant. More fundamental is the existence of the functor $\Theta$ which 
allows to see the various dg-schemes constructed in \cite{ka2,ck1,ck2} as objects in $\mathrm{Ho}(D-Aff^{\sim})$. 

\begin{rmk}
\emph{The above construction of $\Theta$ can be extended from dg-schemes to 
(Artin) dg-stacks.}
\end{rmk}

\end{subsection}

\begin{subsection}{The $D$-stack of $G$-torsors}

As our last example, we present the $D$-stack of $G$-torsors where $G$ is a linear algebraic group $G$. 
As an object in $\mathrm{Ho}(D-Aff^{\sim})$ it is simply $iBG$ (where $BG$ is the usual 
stack of $G$-torsors), but we would like to describe explicitly the
functor $CDGA \longrightarrow SSet$ it represents. \\

Let $H:=\mathcal{O}(G)$ be the Hopf algebra associated to $G$. By considering it 
as an object in the model category of commutative differential graded Hopf algebras, we
can take a cofibrant model $QH$ of $H$, as a dg-Hopf algebra. It is not very hard to check that 
$QH$ is also a cofibrant model for $H$ is the model category of cdga's.
Using the co-algebra structure on $QH$, one sees that the simplicial presheaf
$$Spec\, QH : D-Aff^{op} \longrightarrow SSet$$
has a natural structure of group-like object. In other words, 
$Spec\, QH$ is a presheaf of simplicial groups on $D-Aff$. As the underlying simplicial presheaf of
$Spec\, QH$ is naturally equivalent to $\mathbb{R}Spec\, H\simeq iG$, we will simply denote
this presheaf of simplicial groups by $iG$. 

Next, we consider the category $iG-Mod$, of objects in $D-Aff^{\sim}$ together with an 
action of $iG$. If one sees $iG$ as a monoid in $D-Aff^{\sim}$, the category $iG-Mod$
is simply the category of modules over $iG$. The category $iG-Mod$ is equipped with a notion
of weak equivalences, that are defined through the 
forgetful functor $iG-Mod \longrightarrow D-Aff^{\sim}$ (therefore a morphism 
of $iG$-modules is a weak equivalence if the morphism induced on the underlying objects
is a weak equivalence in $D-Aff^{\sim}$). More generally, there is a model
category structure on $iG-Mod$, such that fibrations and equivalences are defined
on the underlying objects. For any object $F \in iG-Mod$, we also 
get an induced model structure on the comma category $iG-Mod/F$. In particular, it makes
sense to say that two objects $G \longrightarrow F$ and $G' \longrightarrow F$
in $iG-Mod$ are equivalent over $F$, if the corresponding objects
in $\mathrm{Ho}(iG-Mod/F)$ are isomorphic. 

Let $Q$ be a cofibrant replacement functor in the model category $CDGA$.
For any cdga $A$, we have $Spec\, QA \in D-Aff^{\sim}$, the representable $D$-stack represented
by $A \in D-Aff$, that we will consider as $iG$-module for the trivial action.
A $G$\textit{-torsor over} $A$ is defined to be a
$iG$-module $F \in iG-Mod$, together with a fibration of $iG$-modules
$F \longrightarrow Spec\, QA$, such that there exists an \'etale 
covering $A \longrightarrow B$ with the  property that the 
object 
$$F\times_{Spec\, QA}Spec\, QB \longrightarrow Spec\, QB$$ 
is equivalent over $Spec\, QB$ to $iG\times Spec\, QB \longrightarrow Spec\, QB$
(where $iG$ acts on itself by left translations). 

For a cdga $A$, $G$-torsors over $A$ form a full sub-category of $iG-Mod/Spec\, QA$, that will be
denoted by $G-Tors(A)$. This category has an obvious induced notion of weak equivalences, and these equivalences
form a subcategory denoted by $wG-Tors(A)$. Transition morphisms $wG-Tors(A) \longrightarrow
wG-Tors(B)$
can be defined for any morphism $A \longrightarrow B$ by sending
a $G$-torsor $F \longrightarrow Spec\, QA$ to the pull-back
$F\times_{Spec\, QA}Spec\, QB \longrightarrow Spec\, QB$. With a bit of
care, one can make this construction into a (strict) functor 
$$\begin{array}{ccc}
CDGA & \longrightarrow & Cat \\
A & \mapsto & wG-Tors(A).
\end{array}$$

We are now ready to define our functor
$$\begin{array}{cccc}
\mathbb{R}BG : & CDGA & \longrightarrow & SSet \\
& A & \mapsto & |wG-Tors(A)|,
\end{array}$$
where $|wG-Tors(A)|$ is the nerve of the category $wG-Tors(A)$.
The following result says that $\mathbb{R}BG$ is the \textit{associated $D$-stack
to $iBG$} (recall that $BG$ is the Artin stack of $G$-torsors, and that $iBG$ is
its associated $D$-stack defined through the embedding $i$ of Lemma \ref{l2}). 

\begin{prop}\label{p2}
\begin{enumerate}
\item The object $\mathbb{R}BG \in D-Aff^{\sim}$ is a $D$-stack. 

\item There exists an isomorphism $iBG\simeq \mathbb{R}BG$ in the homotopy category $\mathrm{Ho}(D-Aff^{\sim})$.

\end{enumerate}
\end{prop}

An important case is $G=Gl_{n}$, for which we get that the image under $i$ of the stack
$\underline{Vect}_{n}$ of vector bundles of rank $n$ is equivalent to $\mathbb{R}BGl_{n}$ as defined
above. 

\end{subsection}

\end{section}

\begin{section}{The geometry of $D$-stacks}

We are now ready to start our geometric study of $D$-stacks. We will define in this
Section a notion of ($1$)-\textit{geometric} $D$-stack, analogous to the notion
of algebraic stack (in the sense of Artin).
We will also present the theory of \textit{tangent} $D$\textit{-stacks}, as well
as its relations to the cotangent complex. 

\begin{subsection}{Geometricity}

A $1$-geometric $D$-stack is a \textit{quotient of a disjoint union of representable $D$-stacks
by the action of a smooth affine groupoid}. In order to define precisely this notion, we need some preliminaries. 

\begin{enumerate}
\item
Let $f : F \longrightarrow F'$ be a morphism in $\mathrm{Ho}(D-Aff^{\sim})$. We say that $f$ is 
a \textit{representable morphism}, if for any cdga $A$, and any morphism $\mathbb{R}Spec\, A \longrightarrow F'$, 
the homotopy pull-back $F\times_{F'}^{h}\mathbb{R}Spec\, A$ is a  representable $D$-stack (see 
Definition \ref{d6}). 
\item
We say that a $D$-stack $F$ has a representable diagonal if the diagonal morphism 
$\Delta : F \longrightarrow F\times F$ is representable. Equivalently, $F$ has a representable diagonal
if any morphism $\mathbb{R}Spec\, A \longrightarrow F$ from a representable $D$-stack is a representable morphism.
\item
Let $u : A \longrightarrow B$ be a morphism of cdga's. We say that $u$ is \textit{strongly smooth}\footnote{The 
expression \textit{smooth morphism} will be used for a weaker notion in \S 4.4.} if there exists an 
\'etale covering $B \longrightarrow B'$, and a factorization
$$\xymatrix{A \ar[r] \ar[d] & B \ar[d] \\
 A\otimes \mathbb{C}[X_{1},\dots,X_{n}] \ar[r] & B'}$$ 
with $A\otimes \mathbb{C}[X_{1},\dots,X_{n}] \longrightarrow B'$ formally \'etale; 
here $\mathbb{C}[X_{1},\dots,X_{n}]$ is the usual polynomial ring, viewed as a cdga concentrated in degree zero. 
This is an extension of one of the many equivalent characterizations of smoothness for morphisms of schemes 
(see \cite[Prop. 3.24 (b)]{mil}); we learn it from \cite{mcm} 
in which smooth morphisms (called there \textit{thh-smooth}) between $\mathbb{S}$-algebras are defined.
\item A representable morphism of $D$-stacks $f : F \longrightarrow F'$ is called 
\textit{strongly smooth}, if for any morphism from a representable $D$-stack $\mathbb{R}Spec\, A \longrightarrow F'$, the
induced morphism
$$F\times_{F'}^{h}\mathbb{R}Spec\, A \longrightarrow \mathbb{R}Spec\, A$$
is induced by a strongly smooth morphism of cdga's.

\item A morphism $f : F \longrightarrow F'$ in $\mathrm{Ho}(D-Aff^{\sim})$ 
is called a \textit{covering} (or an epimorphism), if the induced morphism $\pi_{0}(F) \longrightarrow \pi_{0}(F')$
is an epimorphism of sheaves. 

\end{enumerate}

Note that definition $(4)$ above makes sense because of $(1)$ and because the functor $A \mapsto \mathbb{R}Spec\, A$
is fully faithful on the homotopy categories. \\

Using these notions, we give the following 

\begin{df}\label{d7}
A $D$-stack $F$ is \emph{strongly ($1$)-geometric} if it satisfies the following two conditions
\begin{enumerate}
\item $F$ has a representable diagonal.

\item There exist representable $D$-stacks $\mathbb{R}Spec\, A_{i}$, and 
a covering 
$$\coprod_{i}\mathbb{R}Spec\, A_{i} \longrightarrow F,$$
such that each of the morphisms $\mathbb{R}Spec\, A_{i} \longrightarrow F$ (which is representable by 1.) is
strongly smooth. Such a family of morphisms will be called a \emph{strongly smooth atlas of $F$}.

\end{enumerate}
\end{df}

\begin{rmk}
\emph{Objects satisfying Definition \ref{d7} are called strongly $1$-geometric $D$-stacks as there 
exists a more general notion of strongly $n$-geometric $D$-stacks, obtained by induction
as suggested in \cite{s}. The notion of strongly $1$-geometric $D$-stacks will be enough
for our purposes (except for our last example in section $5$), 
and we will simply use the expression \textit{strongly geometric $D$-stacks}.}
\end{rmk}

The following proposition collects some of the basic properties of strongly geometric $D$-stacks.

\begin{prop}\label{p3}
\begin{enumerate}
\item Representable $D$-stacks are strongly geometric.

\item The homotopy pull-back of a diagram of strongly geometric $D$-stacks is again
a strongly geometric $D$-stack. In particular strongly geometric $D$-stacks are stable
by finite homotopy limits.

\item If $F$ is any algebraic stack (in the sense
of Artin, see \cite{lm}) with an affine diagonal, then $iF$ is
a strongly geometric $D$-stack.

\item If $F$ is a strongly geometric $D$-stack then $h^{0}(F)$ is an algebraic
stack (in the sense of Artin) with affine diagonal. In particular, 
$ih^{0}(F)$ is again a strongly geometric $D$-stack.

\item For any dg-scheme $(X,\mathcal{A}_{X})$, ($X$ separated and quasi-compact), 
$\Theta(X,\mathcal{A}_{X})$  (see \S 3.3) is a strongly geometric $D$-stack.

\end{enumerate}
\end{prop}

In particular, Proposition \ref{p2} and point $(3)$ above, tell us that 
the derived stack $\mathbb{R}BG$ of $G$-torsors is a 
stongly geometric $D$-stack for any linear algebraic group $G$. \\

We are not going to present the theory in details in this work, but we would like to mention
that standard notions in algebraic geometry (e.g. smooth or flat morphisms, 
sheaves, cohomology \dots) can be extended to strongly geometric $D$
stacks. We refer to \cite{lm} and \cite{s} for the main outline of the constructions. 
The reader will find all details in \cite{hag2}.

\end{subsection}

\begin{subsection}{Modules, linear $D$-stacks and $K$-theory}

Let $\mathbb{G}_{a}$ be the additive group scheme (over $\mathbb{C}$) and consider the object 
$i\mathbb{G}_{a} \in \mathrm{Ho}(D-Aff^{\sim})$. It has
a nice model in $D-Aff^{\sim}$ which is $Spec\, \mathbb{C}[T]$ that we will  denote by $\mathcal{O}$
(note that $\mathbb{C}[T]$ as a cdga in degree $0$ is a cofibrant object).
The $D$-stack $\mathcal{O}$ is actually an 
object in commutative $\mathbb{C}$-algebras, explicitly given by
$$\begin{array}{cccc}
\mathcal{O} : & CDGA & \longrightarrow & (\mathbb{C}-Alg)^{\Delta^{op}} \\
 & A & \mapsto & ([n] \mapsto \Gamma_{n}(A)^{0}),
\end{array}$$
where $\Gamma$ is a fibrant resolution functor. The $D$-stack is called
the \textit{structural $D$-stack}. 

Let us now fix a $D$-stack $F$, and consider the comma category $D-Aff^{\sim}/F$ of $D$-stacks
over $F$; this category is again a model category for the obvious model structure. 
We define the \textit{relative structural $D$-stack} by
$$\mathcal{O}_{F}:=\mathcal{O}\times F \longrightarrow F \in D-Aff^{\sim}/F.$$
Since $\mathcal{O}$ is a $\mathbb{C}$-algebra object, 
we deduce immediately that $\mathcal{O}_{F}$ is also a $\mathbb{C}$-algebra object in the
comma model category $D-Aff^{\sim}/F$. 

Then we can consider the category $\mathcal{O}_{F}-Mod$, of objects in 
$\mathcal{O}_{F}$-modules in the category $D-Aff^{\sim}/F$. If one defines
equivalences and fibrations through the forgetful functor $D-Aff^{\sim}/F \longrightarrow D-Aff^{\sim}$, 
the category $\mathcal{O}_{F}-Mod$ becomes a model category. It has moreover a natural tensor product
structure $\otimes_{\mathcal{O}_{F}}$. 
The model
category $\mathcal{O}_{F}-Mod$ is called the \textit{model category of $\mathcal{O}$-modules on $F$}. \\

Let $A$ be a cdga and $M$ be an (unbounded) $A$-dg module. 
We define a $\mathcal{O}_{Spec\, A}$-module $\widetilde{M}$ in the following way.

Let $\Gamma$ be a fibrant resolution functor on the model category $CDGA$. For any cdga $B$, 
and any integer $n$, we define $\widetilde{M}(B)_{n}$  as 
the set of pairs $(u,m)$, where $u$ is a morphism of cdga's $A \longrightarrow \Gamma_{n}(B)$
(i.e. $u \in Spec\, A(B)$), and $m$ is a degree $0$ element in 
$M\otimes_{A}\Gamma_{n}(B)$ (i.e. $m$ is
a morphism of complexes of $\mathbb{C}$-vector spaces
$m : \mathbb{C} \longrightarrow M\otimes_{A}\Gamma_{n}(B)$). 
This gives a simplicial set $[n] \mapsto \widetilde{M}(B)_{n}$, and therefore
defines an object in $D-Aff^{\sim}$
$$\begin{array}{cccc}
\widetilde{M} :  & CDGA & \longrightarrow & SSet \\
 & B & \mapsto & \widetilde{M}(B).
\end{array}$$

Clearly, the projection $(u,m) \mapsto u$ in the notation above induces a morphism
$\widetilde{M} \longrightarrow Spec\, A$. Finally, this object
is endowed in an obvious way with a structure of $\mathcal{O}_{Spec\, A}$-module. 

This construction, $M \mapsto \widetilde{M}$ induces a functor
$$\widetilde{M} : A-Mod \longrightarrow \mathcal{O}_{Spec\, A}-Mod$$
from the category of (unbounded) dg-$A$-modules, to the category of
$\mathcal{O}_{Spec\, A}$-modules. This functor can be derived (by taking first cofibrant replacements
of both $A$ and $M$) to a functor
$$\mathbb{R}\widetilde{M} : \mathrm{Ho}(A-Mod) \longrightarrow \mathrm{Ho}(\mathcal{O}_{\mathbb{R}Spec\, A}-Mod).$$

\begin{lem}
The functor $\mathbb{R}\widetilde{M}$ defined above is fully faithful. 
\end{lem}

\begin{df}\label{dqcoh}
\begin{enumerate}
\item A $\mathcal{O}$-module on a representable $D$-stack $\mathbb{R}Spec\, A$ is called
\emph{pseudo-quasi-coherent} if it is 
equivalent to some $\mathbb{R}\widetilde{M}$ as above. 

\item 
Let $F$ be a $D$-stack, and $\mathcal{M}$ be a $\mathcal{O}$-module. We say that 
$\mathcal{M}$ is \emph{pseudo-quasi-coherent} if for any 
morphism $u : \mathbb{R}Spec\, A \longrightarrow F$, the pull-back
$u^{*}\mathcal{M}$ is a quasi-pseudo-coherent $\mathcal{O}$-module
on $\mathbb{R}Spec\, A$. 

\end{enumerate}
\end{df}

The construction $M \mapsto \widetilde{M}$ described above also has
a dual version, denoted by $M \mapsto Spel(M)$ and defined in a similar way. 

Let $A$ be a cdga and $M$ be an (unbounded) dg-$A$-module. For a cdga $B$ and an integer 
$n$, we define $Spel(M)(B)_{n}$ to be the set of pairs $(u,\alpha)$, where
$u : A \longrightarrow \Gamma_{n}(B)$ is a morphism of cdga, and
$\alpha : M \longrightarrow \Gamma_{n}(B)$
is a morphism of dg-$A$-modules. This defines a $D$-stack $B \mapsto Spel(M)(B)$
which has a natural projection $(u,\alpha) \mapsto u$, to the $Spec\, A$. Once again, $Spel(M)$
comes equipped with a natural structure of $\mathcal{O}_{Spec\, A}$-module.
Also, this $Spel$ construction 
can be derived, to get a functor
$$\mathbb{R}Spel : \mathrm{Ho}(A-Mod)^{op} \longrightarrow \mathrm{Ho}(\mathcal{O}_{\mathbb{R}Spec\, A}-Mod).$$

\begin{lem}
The functor $\mathbb{R}Spel$ defined above is fully faithful. 
\end{lem}

\begin{df}\label{dlin}
\begin{enumerate}
\item A $\mathcal{O}$-module on a representable $D$-stack $\mathbb{R}Spec\, A$ is called
\emph{representable} if it is 
equivalent to some $\mathbb{R}Spel(M)$ as above. 

\item 
Let $F$ be a $D$-stack, and $\mathcal{M}$ be a $\mathcal{O}$-module. We say that 
$\mathcal{M}$ is \emph{representable} or is a \emph{linear $D$-stack over $F$} if for any 
morphism $u : \mathbb{R}Spec\, A \longrightarrow F$, the pull-back
$u^{*}\mathcal{M}$ is a representable $\mathcal{O}$-module
on $\mathbb{R}Spec\, A$.

\item A \emph{perfect $\mathcal{O}$-module} on a $D$-stack $F$ is
a $\mathcal{O}_{F}$-module which is both pseudo-quasi-coherent and
representable.

\end{enumerate}
\end{df}

One can prove that the homotopy category of perfect $\mathcal{O}$-modules on
$\mathbb{R}Spec\, A$ is naturally equivalent to the full sub-category of $\mathrm{Ho}(A-Mod)$
consisting of strongly dualizable modules, or equivalently of 
dg-$A$-modules which are retracts of finite cell modules (in the sense of \cite[\S III.1]{km}). 
In particular, if $A$ is concentrated in degree $0$, then 
the homotopy category of perfect $\mathcal{O}$-modules on $\mathbb{R}Spec\, A$ is naturally
equivalent to the derived category of bounded complexes of finitely generated projective 
$A$-modules. \\

This notion of perfect $\mathcal{O}$-modules can be used in order to define
the $K$\textit{-theory of }$D$\textit{-stacks}.
For any $D$-stack $F$, one can consider the homotopy category of
perfect $\mathcal{O}$-modules on $F$, that we denote by $D_{\mathrm{Perf}}(F)$. This is a triangulated category 
having a natural Waldhausen model $W\mathrm{Perf}(F)$, from which one can define the $K$-theory spectrum
on the $D$-stack $F$, as $K(F):=K(W\mathrm{Perf}(F))$. The tensor product of $\mathcal{O}$-modules
makes $K(F)$ into an $E_{\infty}$-ring spectrum. Of course, when $X$ is a scheme
$K(iX)$ is naturally equivalent to the $K$-theory spectrum of $X$ as defined
in \cite{tt}. 

A related problem is that of defining reasonable \textit{Chow groups}
and \textit{Chow rings} for strongly geometric $D$-stacks, receiving Chern classes from 
the $K$-theory defined above. We are not aware of any such constructions nor we have 
any suggestion on how to approach the question. It seems however that 
an \textit{intersection theory over} $D$-stacks would be a very interesting 
tool, as it might for example give new interpretations (and probably
extensions) of the notion of \textit{virtual fundamental class} defined in \cite{bef}. 
For this case, the idea would be that for any strongly geometric $D$-stack $F$, there exists
a virtual fundamental class in the Chow group of its truncation $h^{0}F$. 
The structural sheaf  of $F$ should give rise, in the usual way, 
to a fundamental class in its Chow group, 
such that integrating against it over the entire $F$ is the same thing as 
integrating on its truncation $h^{0}F$ against the virtual fundamental class.  
However, even if there is not yet a theory of Chow groups for $D$-stacks, 
if one is satisfied with working with $K$-theory instead of Chow groups, 
the obvious class  $1=:[\mathcal{O}_{F}]\in K_{0}(F)$, 
will correspond exactly to the class of the expected virtual structure sheaf.
\end{subsection}

\begin{subsection}{Tangent $D$-stacks}

Let $Spec\, \mathbb{C}[\epsilon]$ the spectrum of the dual numbers, and let us consider
$iSpec\, \mathbb{C}[\epsilon] \in \mathrm{Ho}(D-Aff^{\sim})$. 

\begin{df}\label{d8}
The \emph{tangent $D$-stack} of a $D$-stack $F$ is defined to be
$$\mathbb{R}TF:=\mathbb{R}\mathcal{HOM}(iSpec\, \mathbb{C}[\epsilon],F) \in \mathrm{Ho}(D-Aff^{\sim}).$$
\end{df}

Note that the zero section morphism $Spec\, \mathbb{C} \longrightarrow Spec\, \mathbb{C}[\epsilon]$
and the natural projection $Spec\, \mathbb{C}[\epsilon] \longrightarrow Spec\, \mathbb{C}$ induces
natural morphisms
$$\pi : \mathbb{R}TF \longrightarrow F \qquad e : F \longrightarrow \mathbb{R}TF,$$
where $e$ is a section of $\pi$. 

An important remark is that for any $D$-stack $F$, the truncation
$h^{0}\mathbb{R}TF$ is equivalent to the tangent stack of $h^{0}F$ (in the sense
of \cite[\S 17]{lm}). In other words, one has
$$h^{0}\mathbb{R}TF \simeq T(h^{0}F).$$
In particular, the $D$-stacks $\mathbb{R}TF$ and $iT(h^{0}F)$ have the same classical points. 
However, it is \textit{not} true in general that $iTF\simeq \mathbb{R}T(iF)$ for a stack $F$. 
Even for a scheme $X$, it is not true that $\mathbb{R}T(iX) \simeq iTX$, except when 
$X$ is smooth. 

\begin{df}\label{d9}
If $x : iSpec\, \mathbb{C} \longrightarrow F$ is a point of a $D$-stack $F$, then the
\emph{tangent $D$-stack of $F$ at $x$} is the homotopy fiber of 
$\pi : \mathbb{R}TF \longrightarrow F$ at the point $x$. It is denoted by
$$\mathbb{R}TF_{x}:=\mathbb{R}TF\times_{F}^{h}iSpec\, \mathbb{C} \in \mathrm{Ho}(D-Aff^{\sim}).$$
\end{df}

Let us now suppose that $F$ is a strongly geometric $D$-stack. One can show that 
$\mathbb{R}TF$ is also strongly geometric. In particular, for any point
$x$ in $F(\mathbb{C})$ the $D$-stack $\mathbb{R}TF_{x}$ is strongly geometric.

Actually much more is true. For any strongly geometric $D$-stack $F$, and any point $x$ in $F(\mathbb{C})$, 
the $D$-stack $\mathbb{R}TF_{x}$ is a \textit{linear} $D$\textit{-stack}
 (over $i\mathrm{Spec}\mathbb{C}$) as defined
in \ref{dlin}. Let us recall that this implies the existence of a natural complex $\mathbb{R}\Omega^{1}_{F,x}$ of 
$\mathbb{C}$-vector spaces (well defined up to a quasi-isomorphism and
concentrated in degree $]-\infty,1]$), with the property that, for any cdga $A$, there exists a natural
equivalence 
$$\mathbb{R}TF_{x}(A)\simeq \mathbb{R}\underline{Hom}_{C(\mathbb{C})}(\mathbb{R}\Omega^{1}_{F,x},A),$$
where $\mathbb{R}\underline{Hom}_{C(\mathbb{C})}$ denotes the mapping space in the model category of
(unbounded) complexes of $\mathbb{C}$-vector spaces. Symbolically, one writes
$$\mathbb{R}TF_{x}=(\mathbb{R}\Omega^{1}_{F,x})^{*},$$
where $(\mathbb{R}\Omega^{1}_{F,x})^{*}$ is the dual complex to $\mathbb{R}\Omega^{1}_{F,x}$. 
In other words, the tangent $D$-stack of $F$ at $x$ ``is'' the
complex $(\mathbb{R}\Omega^{1}_{F,x})^{*}$, which is now concentrated 
in degree $[-1,\infty[$. 

\begin{df}\label{d10}
If $x : iSpec\, \mathbb{C} \longrightarrow F$ is a point of a strongly geometric $D$-stack, then we say that
the \emph{dimension of $F$ at $x$ is defined} if the complex $\mathbb{R}\Omega^{1}_{F,x}$
has bounded and finite dimensional cohomology. If this is the case, the \emph{dimension
of $F$ at $x$} is defined by
$$\mathbb{R}Dim_{x}F:=\sum_{i}(-1)^{i}H^{i}(\mathbb{R}\Omega^{1}_{F,x}).$$
\end{df}

This linear description of $\mathbb{R}TF_{x}$ has actually a global version. In fact, one can define
a \textit{cotangent complex} $\mathbb{R}\Omega^{1}_{F}$ of a strongly geometric $D$-stack, which 
is in general an $\mathcal{O}$-module on $F$ in the sense
of Definition \ref{dqcoh}, which is most of the times quasi-coherent. One then shows that
there exists an equivalence of $D$-stacks over $F$
$$\mathbb{R}TF\simeq \mathbb{R}Spel\, (\mathbb{R}\Omega^{1}_{F}),$$
and in particular that the $D$-stack $\mathbb{R}TF$ is a \textit{linear stack} over $F$ in the
sense of Definition \ref{dlin}.

An already interesting application of this description, is to the case $F=iX$, for 
$X$ a scheme or even an algebraic stack. Indeed, the cotangent
complex $\mathbb{R}\Omega^{1}_{iX}$ mentioned above is
precisely the cotangent complex $\mathbb{L}_{X}$ of \cite[\S 17]{lm}. 
The equivalence 
$$\mathbb{R}T(iX)\simeq \mathbb{R}Spel\, (\mathbb{R}\Omega^{1}_{X})$$
gives a relation between the \textit{purely algebraic} object $\mathbb{L}_{X}$ and
the \textit{geometric object} $\mathbb{R}T(iX)$. In a sense, the usual geometric intuition
about the tangent space is recovered here, at the price of (and thanks to) enlarging the category of objects under study: 
the cotangent complex of a scheme becomes the 
derived tangent space of the scheme considered as a $D$-stack. 
We like to see this
as a possible answer to the following remark of A. Grothendieck (\cite[p. 4]{gr}):\\

[\dots] \textit{Il est tr\`es probable que cette th\'eorie  pourra s'\'etendre de fa\c{c}on \`a donner
 une correspondance entre complexes de chaines de longeur n, et certaines ``n-cat\'egories''  cofibr\'ees 
sur $\underline{C}$; et il n'est pas exclus que par cette voie on arrivera \'egalement \`a une 
``interpr\'etation g\'eom\'etrique'' du complexe cotangent relatif de Quillen.} \\

\end{subsection}

\begin{subsection}{Smoothness}

To finish this part, we investigate various non-equivalent natural notions of smoothness for geometric
$D$-stacks. \\

\textbf{Strong smoothness.} We have already defined the notion of a 
strongly smooth morphisms of cdga's in \S 4.1. We will therefore say that 
a morphism  
$$F \longrightarrow \mathbb{R}Spec\, B$$
from a geometric $D$-stack $F$ is \textit{strongly smooth} if there is  
a strongly smooth atlas $\coprod \mathbb{R}Spec\, A_{i} \longrightarrow F$ as in Definition \ref{d7}, such that 
all the induced morphisms of cdga's $B \longrightarrow A_{i}$
are strongly smooth morphisms of cdga's (see \S 4.1). More generally, 
a morphism between strongly geometric $D$-stacks, $F \longrightarrow F'$, is called
strongly smooth if for any morphism $\mathbb{R}Spec\, B \longrightarrow F'$ the morphism 
$F\times_{F'}^{h}\mathbb{R}Spec\, B \longrightarrow \mathbb{R}Spec\, B$
is strongly smooth in the sense above. 

Strong smoothness is not very interesting for $D$-stacks, as a strongly geometric $D$-stack
$F$ will be strongly smooth if and only if it is of the form $iF'$, for $F'$ 
a \textit{smooth} algebraic stack. \\

\textbf{Standard smoothness.} A more interesting notion is that 
of \textit{standard smooth morphisms}, or simply \textit{smooth morphisms}. 
On the level of cdga's they are defined as follows. 

A morphism of cdga's $A \longrightarrow B$ is called \textit{standard smooth}
(or simply \textit{smooth}), if there exists an \'etale covering
$B \longrightarrow B'$, and a factorization
$$\xymatrix{A \ar[r] \ar[d] & B \ar[d] \\
 A' \ar[r] & B,}$$
such that the $A$-algebra $A'$
is equivalent $A\otimes L(E)$, where $L(E)$ is the free cdga over some bounded complexe
of finite dimensional $\mathbb{C}$-vector spaces $E$. This notion, defined
on cdga's, can be extended (as we did above for strongly smooth morphisms) to 
morphisms between strongly geometric $D$-stacks. 

This notion is more interesting than strong smoothness, as
a strongly geometric $D$-stack can be smooth without being an algebraic stack. However, one
can check that if $F$ is a smooth strongly geometric $D$-stack in this sense, then 
$h^{0}(F)$ is also a smooth algebraic stack. In particular, the derived version of
the stack of vector bundles on a smooth projective surface, discussed in the Introduction (see also conjecture \ref{RVect}),
will never be smooth in this sense as its truncation is the stack of
vector bundles on the surface which is singular in general). 

Nevertheless, smooth morphisms can be used in order to define the following more general 
notion of geometric $D$-stacks.

\begin{df}\label{d7'}
A $D$-stack $F$ is \emph{($1$)-geometric} if it satisfies the following two conditions
\begin{enumerate}
\item The $D$-stack has a representable diagonal.

\item There exists representable $D$-stacks $\mathbb{R}Spec\, A_{i}$, and 
a covering 
$$\coprod_{i}\mathbb{R}Spec\, A_{i} \longrightarrow F,$$
such that each of the morphisms $\mathbb{R}Spec\, A_{i} \longrightarrow F$ is
smooth. Such a family of morphisms will be called a \emph{smooth atlas of $F$}.

\end{enumerate}
\end{df}

Essentially all we have said about strongly geometric $D$-stacks is also
valid for geometric $D$-stacks in the above sense. The typical example of 
a geometric $D$-stack which is not strongly geometric is $BG$, where $G$ is 
a representable group $D$-stack which is not a scheme. For example, 
one can take $G$ to be of the from $\mathbb{R}Spel(M)$ for a non-positively graded 
bounded complex of finite dimensional vector spaces. Then, $G$ is a representable
$D$-stack (it is precisely $\mathbb{R}Spec\, L(M)$, where $L(M)$ is the free
cdga on $M$), and $BG$ is naturally equivalent to $\mathbb{R}Spel(M[-1])$. 
When $M[-1]$ has non-zero $H^{1}$ then $BG$ is not representable anymore but is
$1$-geometric for the above definition. 

More generally, the definition above allows one to consider quotient $D$-stacks
$[X/G]$, where $X$ is a representable $D$-stack and $G$ is a smooth 
representable group $D$-stack acting on $X$. \\

\textbf{fp-smoothness.} The third notion of smoothness is called \textit{fp-smoothness} and
is the weakest of the three and it seems this is the one which is closer to
the smoothness notion referred to in the Derived Deformation Theory program in general. 
It is also well suited in order for the derived stack of vector bundles to be smooth. 

Recall that a morphism of cdga's, $A \longrightarrow B$ is 
\textit{finitely presented} if it is equivalent to a retract of a 
finite cell $A$-algebra, or equivalently if the mapping space
$\mathrm{Map}_{A/CDGA}(B,-)$ commutes with filtered colimits (this is the same
as saying that $\mathbb{R}Spec\, A$ commutes with filtered colimits). 
We will then say that
a morphism of geometric $D$-stacks, $F \longrightarrow F'$ is 
\textit{locally finitely presented} if for any morphism $\mathbb{R}Spec\, A \longrightarrow F'$
there exists a smooth atlas 
$$\coprod \mathbb{R}Spec\, A_{i} \longrightarrow F\times_{F'}\mathbb{R}Spec\, A$$ 
such that 
all the induced moprhisms of cdga's $A \longrightarrow A_{i}$
are finitely presented. Locally finitely presented morphisms will also be
called \textit{fp-smooth} morphisms. The reason for this name is given by the following 
observation.

\begin{prop}\label{p4}
Le $F$ be a geometric $D$-stack which is fp-smooth (i.e. $F \longrightarrow *=i\mathrm{Spec}\mathbb{C}$
is fp-smooth). Then the cotangent complex $\mathbb{R}\Omega^{1}_{F}$ is
a perfect complex of $\mathcal{O}$-modules on $F$. 

In particular, for any point
$x \in F(\mathbb{C})$, the dimension of $F$ at $x$ is defined and locally constant for the
\'etale topology.
\end{prop}

Of course, one has \textit{strongly smooth $\Rightarrow$ smooth $\Rightarrow$
fp-smooth}, but each of these implications is strict. For example, 
a smooth scheme is strongly smooth. Let $E$ be a complex in non-positive degrees which is
cohomologically bounded and of finite dimension. Then $\mathbb{R}Spel(E)$ 
is smooth but not strongly smooth as it is not a scheme in general. Finally, any 
scheme which is a local complete interesection is fp-smooth, but not smooth in general.

\end{subsection}

\end{section}

\begin{section}{Further examples}

In this Section we present three examples of geometric
$D$-stacks: the derived stack of \textit{local systems on a space}, the
derived stack of \textit{vector bundles} and the derived stack
of \textit{associative algebra} and \textit{$A_{\infty}$-categorical structures}.
The derived moduli space of local systems on a space has already been
introduced and defined in \cite{ka2} as a dg-scheme. In the same way, 
the derived moduli space of (commutative) algebra structures 
has been constructed in \cite{ck2} also as a dg-scheme. Finally, 
the formal derived moduli spaces of local systems on a space and
of $A_{\infty}$-categorical structures have been considered 
in \cite{hin2,ko2,kos}. 

The new mathematical content of this part
is the following. First of all we give a construction of the derived moduli stack
of vector bundles, that seems to be new, and we also define
global versions of the formal moduli spaces of $A_{\infty}$-categorical
structures that were apparently not known. We also provide explicit modular descriptions, 
by defining various derived moduli functors, which 
were not known (and probably not easily available), for the constructions of \cite{ka2,ck1,ck2}.

\begin{subsection}{Local systems on a topological space}

Throughout this subsection, $X$ will be a CW-complex. 
For any cdga $A$, we denote by $A-Mod_{X}$ the 
category of presheaves of dg-$A$-modules 
over $X$. We say that a map $\mathcal{M}\rightarrow \mathcal{N}$ 
in $A-Mod_{X}$ is a \textit{quasi-isomorphism} if it induces a 
quasi-isomorphism of dg-$A$-modules on each stalk. This gives a notion of
\textit{equivalences} in the category $A-Mod_{X}$, and of \textit{equivalent objects} 
(i.e. objects which are isomorphic in the localization of the category with respect to equivalences).

A presheaf $\mathcal{M}$ of dg-$A$-modules on $X$ 
will be said \textit{locally on} $X\times A_{\textrm{\'et}}$ \textit{equivalent 
to} $A^{n}$ if, for any $x\in X$, there exists an open neighborhood 
$U$ of $x$ in $X$ and an \'etale cover $A\rightarrow B$, such that 
the presheaves of dg-$B$-modules $\mathcal{M}_{|U}\otimes_{A} B$ and $B^{n}$
are equivalent. We will also say that 
a presheaf $\mathcal{M}$ of dg-$A$-modules is \textit{flat}, if
for any open $U$ of $X$, the dg-$A$-module $\mathcal{M}(U)$ is
cofibrant. By composing with a cofibrant replacement functor in 
$A-Mod$, one can associate to any dg-$A$-module an equivalent flat dg-$A$-module 
(since equivalences are stable by filtered colimits).
The category $w\mathrm{Loc}_{n}(X;A)$ \textit{of 
rank} $n$ \textit{local systems of dg-}$A$\textit{-modules}
has objects those presheaves of flat 
dg-$A$-modules on $X$ which are locally on $X\times 
A_{\textrm{\'et}}$ equivalent to $A^{n}$, and morphisms 
quasi-isomorphisms between them. For morphisms of cdga's
$A \longrightarrow B$ we obtain pull-back functors
$$\begin{array}{ccc} 
w\mathrm{Loc}_{n}(X;A) & \longrightarrow & w\mathrm{Loc}_{n}(X;B) \\
\mathcal{M} & \mapsto & \mathcal{M}\otimes_{A}B.
\end{array}$$
This makes $wLoc_{n}(X;A)$ into a lax functor
from $CDGA$ to categories, that we turn into a strict functor by applying the standard
strictification procedure. 

We denote by $\mathbb{R}\underline{Loc}_{n}(X)$ 
the simplicial presheaf on $D-Aff$ sending a cdga $A$ to 
$\left|w\mathrm{Loc}_{n}(X)\right|$ (the nerve of $w\mathrm{Loc}_{n}(X;A)$). 
We call it the $D$\emph{-pre-stack of rank $n$ derived local systems on} $X$.
 
Obviously, the objects in $w\mathrm{Loc}_{n}(X;A)$ 
are a derived version of the usual local systems of $R$-modules 
on $X$, where $R$ is a commutative ring. More precisely, if we 
consider such an $R$ as a cdga concentrated in degree zero, 
then $R\mathrm{Loc}_{n}(X;R)$ is the closure under quasi-isomorphisms 
of the groupoid of rank $n$ local systems of $R$-modules on $X$; 
in other words, if we invert quasi-isomorphisms in 
the category $w\mathrm{Loc}_{n}(X;R)$ 
then we obtain a category 
which is equivalent to the 
groupoid of rank $n$ local systems of $R$-modules on $X$.

\begin{thm}\label{RLoc}
\begin{enumerate}
	\item The $D$-pre-stack $\mathbb{R}\underline{Loc}_{n}(X)$ 
	is a $D$-stack. Furthermore, one has
	$\mathbb{R}\underline{Loc}_{n}(pt) \simeq iBGl_{n}$.
	\item One has an equivalence
	$$h^{0}\mathbb{R}\underline{Loc}_{n}(X)\simeq [Hom(\pi_{1}(X),Gl_{n})/Gl_{n}],$$ 
	between the truncation of $\mathbb{R}\underline{Loc}_{n}(X)$ and
	the (Artin) stack of local systems on $X$.  
	\item If $S(X)$ denotes the singular complex of $X$, 
	we have the following isomorphisms in $\mathrm{Ho}(D-Aff^{\sim})$, 
	$$\mathbb{R}\underline{Loc}_{n}(X) \simeq 
	\mathbb{R}\mathcal{HOM}(\underline{S(X)}, iBGl_{n})\simeq  
	\mathbb{R}\mathcal{HOM}(\underline{S(X)}, \mathrm{Loc}_{n}(\mathrm{pt})),$$
where $\mathbb{R}\mathcal{HOM}$ denotes the Hom-stack (internal 
Hom in $\mathrm{Ho}(D-Aff^{\sim})$) and $\underline{S(X)}$ denotes the simplicial 
constant presheaf with value $S(X)$. 
	\item For any rank $n$ local system $L$ on $X$, the tangent $D$-stack of
	$\mathbb{R}\underline{Loc}_{n}(X)$ at $L$ is the complex
	$C^{*}(X,\underline{End}(L))[1]$, of cohomology of $X$ with coefficients in $\underline{End}(L)$.
	\item If $X$ is a \emph{finite} CW-complex, then the stack  
	$\mathbb{R}\underline{Loc}_{n}(X)$ is strongly geometric, 
	fp-smooth  of (the expected) dimension $-n^{2}\chi(X)$, 
	$\chi(X)$ being the Euler characteristic of $X$.
\end{enumerate}
\end{thm}

Note that the classical points of $\mathbb{R}\underline{Loc}_{n}(X)$ (i.e. morphisms from $i\mathrm{Spec}k$, 
for some commutative ring $k$) coincide with the classical points of its 
truncation $h^{0}\mathbb{R}\underline{Loc}_{n}(X)$ which coincides 
with the usual (i.e. not derived) stack of rank $n$ local systems on $X$. So 
we have no new classical points, as desired.

Let us give only some remarks to show what the proof of Theorem \ref{RLoc} 
really boils down to. First of all notice that the first assertion is a 
consequence of the second one, once one knows that $\mathbb{R}\underline{Loc}_{n}(\mathrm{pt})\simeq iBGL_{n}$ and 
is a stack; so we are reduced to prove the absolute case ($X=\mathrm{pt}$) of 1. and 2. 
The first two properties in 3. follows from 2., the finiteness of $X$ and the analogous 
properties of $BGl_{n}$. Finally the dimension count in 3. is made by an explicit 
computaion of the tangent $D$-stack 
at some local system $E$. Explicitly, 
one finds that (in the notations of $\S 4.3$) 
$(\mathbb{R}\Omega^{1}_{\mathbb{R}\underline{Loc}_{n}(X),E})^{*}$ is the 
complex $C^{*}(X,\underline{End}(E))[1]$, which is a complex of $\mathbb{C}$-vector 
spaces concentrated in degrees $[-1,\infty[$ whose Euler characteristic is 
exactly $-n^{2}\chi(X)$.

\begin{rmk}
\emph{ The example of local systems is one of those cases where there is a 
canonical way to \textit{derive} the usual moduli stack (see the discussion in Section 6). In fact, in this 
case we have 
$\mathcal{HOM}(\underline{S(X)},\mathrm{Loc_{n}}(\mathrm{pt}))
\simeq \mathrm{Loc}_{n}(X)$, for any CW-complex $X$, where $\mathcal{HOM}$ 
denotes the (underived) Hom-stack between (underived) stacks; therefore the natural thing to do is 
to first view the usual absolute stack $\mathrm{Loc}_{n}(pt)$ as a derived stack 
via the inclusion $i$ and then derive the Hom-stack from $\underline{S(X)}$ to 
$i\mathrm{Loc}_{n}$. This authomatically gives an extension of 
$\mathrm{Loc}_{n}(X)$ i.e. a canonical derivation of it.}
\end{rmk} 

It is important to notice that the $D$-stack $\mathbb{R}\underline{Loc}_{n}(X)$ 
might be non-trivial even if $X$ is simply connected. Indeed, the tangent at the
unit local system is always the complex $C^{*}(X,\mathbb{C})[1]$. 
This shows that $\mathbb{R}\underline{Loc}_{n}(X)$ contains interesting information
concerning the higher homotopy type of $X$. As noticed in the Introduction of \cite{kps}, this is one of the reasons why
the $D$-stack $\mathbb{R}\underline{Loc}_{n}(X)$ might be an interesting object in order
to develop a version of non-abelian Hodge theory. We  will therefore 
ask the same question as in \cite{kps}. 

\begin{Q}
Let $X$ be a smooth projective complex variety and $X^{top}$ its underlying topological
space. Can one extend the non-abelian Hodge structure defined on the moduli space
of local systems in \cite{s2}, to some kind of Hodge structure on the whole
$\mathbb{R}\underline{Loc}_{n}(X)$ ?
\end{Q}

This question is of course somewhat imprecise, and it is not clear that the object $\mathbb{R}\underline{Loc}_{n}(X)$ 
itself could really support an interesting Hodge structure.
However, we understand the previous question in a much broader sense, as
for example it includes the question of defining derived versions of the
moduli spaces of flat and Higgs bundles, and to study their 
relations from a non-abelian Hodge theoretic point of view, as done for example in \cite{s2}.

\end{subsection}

\begin{subsection}{Vector bundles on a projective variety}

We now turn to the example of the derived stack of vector bundles, which is
very close to the previous one. Let $X$ be a fixed 
smooth projective variety. 

If $A$ is a cdga, we consider the space $X$ (with the Zariski topology) together with its
presheaf of cdga $\mathcal{O}_{X}\otimes A$. It makes sense to consider
also presheaves of dg-$\mathcal{O}_{X}\otimes A$-modules on $X$ and morphisms
between them. We define 
a notion of equivalences between such presheaves, by saying the 
$f : \mathcal{M} \longrightarrow \mathcal{N}$ is an equivalence
if it induces a quasi-isomorphism at each stalks. Using this notion of equivalences
we can talk about equivalent dg-$\mathcal{O}_{X}\otimes A$-modules (i.e. objects which become 
isomorphic in the localization of the category with respect to quasi-isomorphisms).

We say that a presheaf of dg-$\mathcal{O}_{X}\otimes A$-module $\mathcal{M}$ on $X$
is a \textit{vector bundle of rank} $n$, if locally on $X_{zar}\times A_{\textrm{\'et}}$ it is equivalent to 
$(\mathcal{O}_{X}\otimes A)^{n}$ (see the previous Subsection for details on this definition). We consider
the category $wVect_{n}(X,A)$, of dg-$\mathcal{O}_{X}\otimes A$-modules which are
vector bundles of rank $n$ and flat (i.e. for each 
open $U$ in $X$, the $\mathcal{O}_{X}(U)\otimes A$-module $\mathcal{M}(U)$ is
cofibrant), and equivalences between them. By the standard strictification procedure we obtain 
a presheaf of categories
$$\begin{array}{ccc}
CDGA & \longrightarrow & Cat \\
A & \mapsto & wVect_{n}(X,A) \\
(A \rightarrow B) & \mapsto & (\mathcal{M} \mapsto \mathcal{M}\otimes_{A}B).
\end{array}$$
We then deduce a simplicial presheaf by appying the nerve construction
$$\begin{array}{cccc}
\mathbb{R}\underline{Vect}_{n}(X) : & CDGA & \longrightarrow & Cat \\
 & A & \mapsto & |wVect_{n}(X,A)|. 
\end{array}$$
This gives an object $\mathbb{R}\underline{Vect}_{n}(X) \in D-Aff^{\sim}$ that we call
the derived moduli stack of rank $n$ vector bundles on $X$. \\

We state the following result as a conjecture, as we have not checked all details. 
However, we are very optimistic about it, as we think that a proof will probably consist
of reinterpreting the constructions of \cite{ck1} in our context.

\begin{conj}\label{RVect}
\begin{enumerate}
	\item The $D$-pre-stack $\mathbb{R}\underline{Vect}_{n}(X)$ 
	is a strongly geometric, fp-smooth $D$-stack. 
	
	\item  There exists a natural isomorphism in $\mathrm{Ho}(D-Aff^{\sim})$
	$$\mathbb{R}\underline{Vect}_{n}(X) \simeq \mathbb{R}\mathcal{HOM}(X,iBGl_{n}).$$
	 
	\item One has an equivalence
	$$h^{0}\mathbb{R}\underline{Vect}_{n}(X)\simeq \underline{Vect}_{n}(X)$$
	between the truncation of the $D$-stack $h^{0}\mathbb{R}\underline{Vect}_{n}(X)$
	and the (Artin) stack of rank $n$ vector bundles on $X$.
	 
	\item The tangent $D$-stack of $\mathbb{R}\underline{Vect}_{n}(X)$ at a
	vector bundle $E$ on $X$, is the complex
	$$C^{*}(X_{Zar},\underline{End}(E))[1].$$
\end{enumerate}
\end{conj}

The same remark as in the case of the derived stack of local systems
holds. Indeed, the usual Artin stack of vector bundles on $X$ is given by
$\mathbb{R}\mathcal{HOM}(X,BGl_{n})$, and our $D$-stack 
of vector bundles on $X$ is $\mathbb{R}\mathcal{HOM}(iX,iBGl_{n})$.

\end{subsection}

\begin{subsection}{Algebras and $A_{\infty}$-categorical structures}

In this last Subsection we present the derived moduli stack of associative algebra structures and
$A_{\infty}$-categorical structures. These are global versions of the formal moduli
spaces studied in \cite{ko2,kos}. \\

\textbf{Associative algebra structures.} We  are going to 
construct a $D$-stack $\mathbb{R}Ass$, classifying associative $dg$-algebra structures.

Let $A$ be any cdga, and let us consider the category of (unbounded) 
associative differential graded $A$-algebras 
$A-Ass$ (i.e. $A-Ass$ is the category of
monoids in the symmetric monoidal category $A-Mod$, of (unbounded) dg-$A$-modules)\footnote{By definition
our associative $A$-dga's are then all central over $A$ since they are 
commutative monoids in $A-Mod$.}. 
This category is a model category for which the weak equivalences
are the quasi-isomorphisms and fibrations are epimorphisms. We restrict ourselves
to the category of cofibrant objects $A-Ass^{c}$, 
and consider the sub-category $wA-Ass^{c}$ consisting of equivalences 
only. If $A \longrightarrow A'$ is any morphism of
cdga's, then we have pull-back functors
$$\xymatrix{wA-Ass^{c} \ar[rr]^-{-\otimes_{A}A'} &  & wA'-Ass^{c}}.$$
This defines a (lax) functor on the category of cdga's that we immediately strictify by
the standard procedure. We will therefore assume that the above constructions are
strictly functorial in $A$. By passing to the corresponding nerves we get 
a presheaf of simplicial sets
$$\begin{array}{cccc}
\mathbb{R}Ass : & CDGA & \longrightarrow & SSet  \\
 & A & \mapsto & |wA-Ass^{c}|. 
\end{array}$$
This gives a well defined object $\mathbb{R}Ass$ in $D-Aff^{\sim}$. 

We define a sub-simplicial presheaf $\mathbb{R}Ass_{n}$ of $\mathbb{R}Ass$, consisting
of associative dg-$A$-algebras $B$ for which there exists an 
\'etale covering $A \longrightarrow A'$ such that the dg-$A'$-module
$B\otimes_{A}^{\mathbb{L}}A'$ is equivalent to $(A')^{n}$.

\begin{thm}
\begin{enumerate}
\item The $D$-pre-stack $\mathbb{R}Ass_{n}$ 
is a $D$-stack.

\item The $D$-stack $\mathbb{R}Ass_{n}$ is strongly geometric. Furthermore, 
$h^{0}\mathbb{R}Ass_{n}$ is naturally equivalent to the (usual) Artin stack
of associative algebra structures on $\mathbb{C}^{n}$.

\item For any global point $V : *  \longrightarrow \mathbb{R}Ass_{n}$, corresponding 
to an associative $\mathbb{C}$-algebra $V$, the tangent $D$-stack
of $\mathbb{R}Ass_{n}$ at $V$ is the complex 
$\mathbb{R}Der(V,V)[1]$ 
of (shifted) derived derivations from $V$ to $V$. 

\end{enumerate}
\end{thm}

From $(3)$ we see that the geometric $D$-stack $\mathbb{R}Ass_{n}$ is not
fp-smooth. Indeed, Quillen gives in \cite[Ex. 11.8]{q} an example of 
a point in $\mathbb{R}Ass_{n}$ at which the dimension in the sense of Definition \ref{d10} is not defined. \\

The previous theorem can also be extended in the following way. Let $V$ be a fixed cohomologically
bounded and finite dimensional complex of $\mathbb{C}$-vector spaces. We define
$\mathbb{R}Ass_{V}$ to be the sub-simplicial presheaf of $\mathbb{R}Ass$ consisting
of associative dg-$A$-algebras $B$ for which there exists 
an \'etale covering $A \longrightarrow A'$ such that the
dg-$A'$-module $B\otimes^{\mathbb{L}}_{A}A'$ is equivalent to 
$A'\otimes V$. 

One can show that $\mathbb{R}Ass_{V}$ is again a $D$-stack, but it is not in general strongly geometric
in the sense of Definition \ref{d7} (nor in the sense
of Definition \ref{d7'}). However, we would like just to mention that 
$\mathbb{R}Ass_{V}$ is still \textit{geometric} in some sense when considered
as a \textit{stack over unbounded cdga's} (the reader will find details in the forthcoming paper \cite{hag2}).
The tangent $D$-stack of $\mathbb{R}Ass_{V}$ at a point is given by the same formula as before.

The construction of $\mathbb{R}Ass_{V}$ can also be extended to classify 
algebra structures over an operad on the complex $V$. One can check that 
the $D$-stacks one obtains in this way are again geometric. These are 
the geometric counterparts of the (discrete) moduli spaces described by C. Rezk in \cite{re}. \\

\textbf{$A_{\infty}$-Categorical structures\footnote{We are working here with the
stronger notion of dg-category (or \textit{strict $A_{\infty}$-categories}), and of course one
could also use $A_{\infty}$-categories instead. However, as the homotopy theories of dg-categories
and of $A_{\infty}$-categories are equivalent, the $D$-stacks obtained would be the same.}.} 
Let $A$ by any cdga.
Recall that a \textit{dg-}$A$\textit{-category} $C$ consists of the following data
\begin{enumerate}
\item A set of objects $Ob(C)$.
\item For each pair of object $(x,y)$ in $Ob(C)$, 
a (unbounded) dg-$A$-module $C_{x,y}$.
\item For each triplet of object $(x,y,z)$ in $Ob(C)$, 
a composition morphism $C_{x,y}\otimes_{A}C_{y,z} \longrightarrow C_{x,z}$
which satisfies obvious associativity and unital conditions. 
\end{enumerate}

\bigskip

There is an obvious notion of \textit{morphism} between dg-$A$-categories. There is also 
a notion of \textit{equivalences} of dg-$A$-categories: they are morphisms
$f : C \longrightarrow C'$ satisfying the following two conditions
\begin{enumerate}
\item For any pair of objects $(x,y)$ of $C$, the induced morphism
$f_{x,y} : C_{x,y} \longrightarrow C'_{x,y}$ is
a quasi-isomorphism of dg-$A$-modules.

\item Let $H^{0}(C)$  (resp. $H^{0}(C')$) be the categories having 
repectively the same set of objects as $C$ (resp. as $C'$), and
$H^{0}(C_{x,y})$ (resp. $H^{0}(C'_{x,y})$) as set of morphisms 
from $x$ to $y$. Then, the induced morphism
$$H^{0}(f) : H^{0}(C) \longrightarrow H^{0}(C')$$
is an equivalence of categories (in the usual sense).

\end{enumerate}

Using these definitions, one has for any cdga $A$, a category 
$A-Cat$ of dg-$A$-categories, with a sub-category of equivalences $wA-Cat$. 
Furthermore, if $A \longrightarrow A'$ is a morphism of cdga, one
has a pull-back functor
$A-Cat \longrightarrow A'-Cat$, obtained
by tensoring the dg-$A$-modules $C_{x,y}$ with $A'$. 
With a bit of care (e.g. by restricting to \textit{cofibrant}
dg-$A$-categories), one gets a simplicial presheaf
$$\begin{array}{cccc}
\mathbb{R}Cat : & CDGA & \longrightarrow & SSet \\
 & A & \mapsto & |wA-Cat|,
\end{array}$$
that is an object in $D-Aff^{\sim}$.

We now fix a \textit{graph} $\mathcal{O}$ of non-positively graded complexes of $\mathbb{C}$-vector spaces. This means that 
$\mathcal{O}$ is the datum of a set $\mathcal{O}_{0}$, and for any $(x,y) \in \mathcal{O}$, of a complex $\mathcal{O}_{x,y}$. We will assume that all the complexes 
$\mathcal{O}_{x,y}$ are bounded with finite dimensional cohomology.
We consider the sub-simplicial presheaf $\mathbb{R}Cat_{\mathcal{O}}$ of $\mathbb{R}Cat$, consisting
of all those dg-$A$-categories $C$ such that locally on $A_{\textrm{\'et}}$ the underlying graph of $C$ 
is equivalent to $\mathcal{O}\otimes A$; the underlying graph of $C$ is defined
to be the graph $G(C)$ whose set of objects is a set of representatives 
of isomorphism classes of objects in $H^{0}(C)$, and whose complexes 
of morphisms are the ones of $C$. The simplicial presheaf $\mathbb{R}Cat_{\mathcal{O}}$ classifies
dg-categorical structures on the graph $\mathcal{O}$. 

The following theorem identifies the tangent of $\mathbb{R}Cat_{\mathcal{O}}$.

\begin{thm}
Let $\widetilde{\mathbb{R}Cat_{\mathcal{O}}}$ be the associated $D$-stack to the
$D$-pre-stack $\mathbb{R}Cat_{\mathcal{O}}$. 
For any global point $C : *  \longrightarrow \mathbb{R}Cat_{\mathcal{O}}$, corresponding 
to a dg-category $C$, the tangent $D$-stack
of $\widetilde{\mathbb{R}Cat_{\mathcal{O}}}$ at $C$ is the whole (shifted) Hochschild cohomology complex 
$C^{*}(C,C)[2]$ \emph{(see e.g. \cite[2.1]{kos} or \cite[2]{so})}.
\end{thm}

\begin{rmk}
\emph{For a cdga $A$, points in $\widetilde{\mathbb{R}Cat_{\mathcal{O}}}(A)$ can be described
as certain \textit{twisted forms} of dg-$A$-categories on the \'etale site of $A$.}
\end{rmk}

Let us suppose that $\mathcal{O}$ is now a graph of finite dimensional vector spaces
(i.e. the complexes $\mathcal{O}_{x,y}$ are concentrated in degree $0$ for any $x,y$). Then one can show that 
the $D$-stack $\widetilde{\mathbb{R}Cat_{\mathcal{O}}}$ is \textit{strongly $2$-geometric}.
Here we use a notion of strongly $n$-geometric $D$-stacks obtained by iterating 
Definition \ref{d7}. The reader will find details about higher geometric stacks
in \cite{hag2} and might also wish to consult \cite{s}.
Note that the $D$-stack $\widetilde{\mathbb{R}Cat_{\mathcal{O}}}$ cannot be $1$-geometric, 
as its truncation
$h^{0}\widetilde{\mathbb{R}Cat_{\mathcal{O}}}$ is the ($2$-)stack
of linear categories. As a $1$-geometric (not derived) stack is always $1$-truncated (as opposite to the
derived case), this shows that $\widetilde{\mathbb{R}Cat_{\mathcal{O}}}$ must be at least
$2$-geometric. 

As in the case of $\mathbb{R}Ass$, if the graph $\mathcal{O}$ is not a graph of vector spaces, then
the $D$-stack $\widetilde{\mathbb{R}Cat_{\mathcal{O}}}$ is not strongly $2$-geometric anymore, but is
still geometric in some sense, when considered as a stack over unbounded cdga's. \\

Let $V$ be a bounded complex with finite dimensional cohomology, also considered as a graph
of complexes with a unique object.
Then, there exists a natural morphism
$$\mathbb{R}Ass_{V} \longrightarrow \widetilde{\mathbb{R}Cat_{V}},$$
that sends an associative dga to the dg-category, with one object, it defines. 
This morphism is actually a \textit{gerbe} in the following sense. 
If $B : \mathbb{R}Spec\, A \longrightarrow \mathbb{R}Ass_{V}$ corresponds to 
an associative $A$-dga $B$, then the homotopy fiber $F$  of the previous morphism is locally equivalent to
the $D$-stack over $\mathbb{R}Spec\, A$ sending a cdga $A \rightarrow A'$ 
to the simplicial set $K((B\otimes_{A}A')^{*},1)$, where
$(B\otimes_{A}A')^{*}$ is the loop space of invertible elements in $B\otimes_{A}A'$
(i.e. the mapping space $\mathrm{Map}_{A'-alg}(A'[T,T^{-1}],B\otimes_{A}A')$). In particular, one
deduces that the morphism $\mathbb{R}Ass_{V} \longrightarrow \widetilde{\mathbb{R}Cat_{V}}$ is 
a \textit{smooth fibration} of $D$-stacks. 
This smooth morphism induces in particular an exact triangle between the
tangent $D$-stacks
$$\xymatrix{\mathbb{R}TF_{B} \ar[r] & \mathbb{R}T(\mathbb{R}Ass_{V})_{B} \ar[r] & 
\mathbb{R}T(\widetilde{\mathbb{R}Cat_{V}})_{B} \ar[r]^-{+1} & }$$
which can also be written
$$\xymatrix{ B[1] \ar[r] & \mathbb{R}Der_{A}(B,B)[1] \ar[r] & 
C^{+}_{A}(B,B)[2] \ar[r]^-{+1} & }$$
which is our way of understanding the triangle appearing in \cite[p. 59]{ko2} (at least
for $d=1$).

The fact that the tangent $D$-stack of $\widetilde{\mathbb{R}Cat_{V}}$ at a dg-category with only one
object is the whole (shifted) Hochschild complex $C^{*}(A,A)[2]$, where $A$ is the dg-algebra
of endomorphisms of the unique object, is also our way to understand 
the following sentence from \cite[p. 266]{kos}.\\

\textit{In some sense, the full Hochschild complex controls
deformations of the $A_{\infty}$-category with one object, such that its endomorphism space is
equal to $A$.} \\

\bigskip

We see that the previous results and descriptions produce 
global versions of the formal moduli spaces of 
$A_{\infty}$-categories studied for example in \cite{kos,so}. This also shows that 
there are interesting higher geometric stacks, and probably 
even more interesting examples will be given by the $D$-stack
of \textit{$(n-1)$-dg-categories} (whatever these are) as suggested by a higher analog of the exact 
triangle above (see \cite[2.7 \, Claim 2]{ko2}).
 
\end{subsection}

\end{section}

\begin{section}{Final comments on deriving moduli functors}

In this Subsection, keeping in mind the examples presented before, we would 
like to discuss, from a more general point of view, the problem of 
derivation of moduli functors, with the aim of at least making explicit 
some general features shared by the examples.\\

Suppose $M:(Aff)^{\mathrm{op}}=(\mathbb{C}-alg)\longrightarrow Set$ is 
a functor arising from some geometric moduli problem e.g., the problem 
of classifying isomorphism classes of families of (pointed) curves of 
a given genus. Very often, the \textit{moduli functor} $M$ is not representable and only
admits a coarse moduli space. As its name says, when passing to a coarse moduli space some
 information is lost. The theory of \textit{stacks in groupoids} was originally invented to correct this annoyance, by 
looking at natural \textit{extensions} of $M$, i.e. to functors $\mathcal{M}_{1}$, from 
algebras to groupoids, such that the following diagram commutes
$$
\xymatrix{Aff^{\mathrm{op}}\ar[r]^{M} \ar[dr]_{\mathcal{M}_{1}} & Set \\
                                                            & Grpd, \ar[u]^-{\pi_{0}} }
$$
Here the 
vertical arrow assigns to a given groupoid its set of isomorphisms classes of objects. 
Of course, the point of the theory of stacks in groupoids is precisely to 
develop a \textit{geometry} on this kind of functors.

More generally, other natural \textit{higher} moduli problems are not representable
even when considered as stacks in groupoids, e.g. the $2$-stack perfect complexes of
length $1$, the $2$-stack of linear categories \dots; the theory of \textit{higher stacks}
precisely says that one should consider $M$ extended as follows
$$
\xymatrix{Aff^{\mathrm{op}}\ar[r]^{M} \ar[rdd]_-{\mathcal{M}} \ar[dr]^{\mathcal{M}_{1}} & Set \\
                                                            & Grpd \ar[u]_-{\pi_{0}} \\
                                                            & SSet, \ar[u]_-{\Pi_{1}}}
$$
where the functor $\Pi_{1}$ maps a simplicial sets to its fundamental groupoid. The notion
of geometric $n$-stacks of \cite{s} can then be used in order to set up a \textit{geometry} over these
kind of objects, in pretty much the same way one is doing geometry over stacks in groupoids. \\
 
The idea of derived algebraic geometry is to seek for derived extensions 
of $M$, $\mathcal{M}_{1}$ and $\mathcal{M}$ i.e. to extend not (only) the target
category of this functors but more crucially the source category in 
a ``derived'' direction.  More precisely, we define a \textit{derived
extension} of a functor $\mathcal{M} : Aff^{op} \longrightarrow SSet$, as above,  
to be a functor 
$\mathbb{R}\mathcal{M}:  (D-Aff)^{\mathrm{op}}\longrightarrow SSets$
making the following diagram commute

$$
\xymatrix{
Aff^{\mathrm{op}}  \ar[dd]_{j} \ar[r]^{M} \ar[dr]^{\mathcal{M}_{1}} \ar[rdd]_-{\mathcal{M}}   & Set \\
                                                                & Grpd \ar[u]_-{\pi_{0}} \\                                                                
          D-Aff^{\mathrm{op}} \ar[r]_{\mathbb{R}\mathcal{M}}  & SSet \ar[u]_-{\Pi_{1}}              }                                                 
$$
where $j$ denotes the natural inclusion (a $\mathbb{C}$-algebra viewed as a 
cdga concentrated in degree zero). The above diagram shows that, 
for any derived extension $\mathbb{R}\mathcal{M}$, we have
$$
\pi_{0}\mathbb{R}\mathcal{M}(j(\mathrm{Spec}R))\simeq M(\mathrm{Spec}R)
$$
and moreover
$$
\Pi_{1}\mathbb{R}\mathcal{M}(j(\mathrm{Spec}R))\simeq \mathcal{M}_{1}(\mathrm{Spec}R)
$$
for any commutative $\mathbb{C}$-algebra $R$. In other words, the 0-truncation of 
$\mathbb{R}\mathcal{M}$ gives back $M$ when restricted to the image of $j$, 
while its 1-truncation gives back $\mathcal{M}_{1}$. \\

What about the existence or uniqueness of a derived extension $\mathbb{R}\mathcal{M}$ ? 
First of all, extensions always exists, as one can use the \textit{trivial} one given by
the functor $i$ of \S 3.2. But of course, this extension is far from being unique
and usually does not give the \textit{expected answer}. However, there is no canonical choice
for an extension which could be \textit{nicer} than others. This tells us that 
the choice of the extended moduli functor $\mathbb{R}\mathcal{M}$ highly depends
on the geometrical meaning of the original moduli functor $M$, $\mathcal{M}_{1}$
of $\mathcal{M}$. We would like to give here a clear example to show this. 

Let $S^{2}$ be the $2$-dimensional
sphere, and let us consider $\mathcal{M}_{1}:=\underline{Loc}_{n}(S^{2})$, the moduli stack of
rank $n$ local systems on $S^{2}$. We clearly have $\mathcal{M}_{1}\simeq BGL_{n}$. 
If one thinks of $\mathcal{M}_{1}$ simply as $BGl_{n}$, and forget about the fact that 
it is the moduli stack of local systems on $S^{2}$, then a reasonable extension of
$\mathcal{M}_{1}$ is simply $iBGl_{n}\simeq \mathbb{R}BGl_{n}$ as described in \S 3.4. However, if
one remembers that $\mathcal{M}_{1}$ is $\underline{Loc}_{n}(S^{2})$, then the \textit{correct}
(or at least \textit{expected}) extension is $\mathbb{R}\underline{Loc}_{n}(S^{2})$ presented
in Theorem \ref{RLoc}. Definitely, these two extensions are very different. This example
shows that the \textit{expected} extension $\mathbb{R}\mathcal{M}$ depends very much on the way we
\textit{think} of the original moduli problem $\mathcal{M}$. In a way we are more
\textit{deriving our interpretation} of the moduli functor rather than the moduli functor itself. 
Another example of the existence of multiple choices can be found in \cite{ck2}, in which 
the derived Hilbert dg-scheme is not the same as the derived $Quot_{\mathcal{O}}$ dg-scheme. \\

Nevertheless, the derived extension of a moduli functor 
that typically occurs in algebraic geometry, is expected to satisfy certain properties and this gives some serious hints in order
to guess the \textit{correct answer}.  

First of all, in general, one knows a priori what is
the \textit{expected derived tangent stack} $T\mathbb{R}\mathcal{M}$ 
(or, at least, the disembodied derived tangent complexes at the points, 
the $(\Omega^{1}_{\mathbb{R}\mathcal{M},x})^{*}$'s in the notations of \S 4.3); 
namely, this is true in the case where $\mathcal{M}$ classifies vector bundles 
over a scheme, local systems over a topological space, families curves 
or higher dimensional algebraic varieties, stable maps from a fixed scheme 
and so on. For some examples of the expected derived tangent spaces 
we refer again to \cite{ck1,ck2}. To put it slightly differently, it is 
always the case that one looks for a derived extension by 
simply requiring it to have the expected derived tangent 
stack. This is essentially due to the fact that the correct derived 
deformation theory of the moduli problem has already been guessed, and the
corresponding, already established, formal theoy is based on this guess
(see \cite{hin2,kos,so}, to quote a few).

Even if this does not say exactly how to construct
a derived extension, it certainly puts some constraints on the possible choices. 
To go a bit further, one may notice that all the 
usual moduli functor occurring in algebraic geometry classify 
\textit{families of geometric objects} over varying base schemes. 
To produce a derived extension $\mathbb{R}\mathcal{M}$, the main principle is then the following \\

\textbf{Main principle:} \textit{Let $\mathcal{M}$ be a moduli stack classifying certain kind of 
families of geometric objects over varying commutative algebras $A$. In order to guess what the extended
moduli stack $\mathbb{R}\mathcal{M}$ should be, guess first what is a family of geometric
objects of the same type parametrized by a commutative dga $A$.} \\

\bigskip

In the case, 
for example, where the classical notion of a family is defined 
through the existence of a map with some properties (like for 
example in the case of the stack of curves), the derived analog 
is more or less clear: one should say that a derived family over a cdga $A$ is 
just a map of simplicial presheaves $F \longrightarrow \mathbb{R}Spec\, A$, having the same properties in the derived 
sense (e.g., as we extended the notion of \'etale morphism of schemes to cdga's, see \S 2.2, 
the same can be done with the notions of smooth, flat \dots morphisms of schemes). Then, a natural 
candidate for a derived notion of family of geometric objects, is given by \textit{any} 
derived analog of a family such that \textit{when restricted along} 
$\mathrm{Spec}(H^{0}(A))\longrightarrow \mathrm{Spec}A$ it becomes \textit{equivalent} 
to an \textit{object coming from} $\mathcal{M}(\mathrm{Spec}(H^{0}(A)))$. This 
condition, required in order to really get a derived extension, essentially says that the derived version of a family of geometric objects 
should reduce to a non-derived family of geometric objects in the non-derived or 
scheme-like direction, i.e. along $\mathrm{Spec}(H^{0}(A))\longrightarrow \mathrm{Spec}A$. A typical example of this case is the one of $G$-torsors given  
in \S 3.4.  Another example would be that of the moduli stack of surfaces. One could say for example that 
a smooth projective family of surfaces over a cdga $A$, is a strongly smooth morphism of strongly geometric $D$-stacks
$F \longrightarrow \mathbb{R}Spec\, A$, such that for any geometric point $x : Spec\, \mathbb{C} \longrightarrow \mathbb{R}Spec\, A$, 
the pull-back $F\times^{h}_{\mathbb{R}Spec\, A}Spec\, \mathbb{C}$ is equivalent to a smooth projective surface. \\

Though this gives perhaps only a vague recipe of a possible construction of derived 
extensions of some of the moduli functors occurring in algebraic geometry, we thought it 
was worthwhile presenting it, if not certainly to solve the problem at least to pose it in a general perspective.

\end{section}

\end{document}